\documentstyle[12pt]{article}
\catcode `\@=11
\@addtoreset{equation}{section}

\def\harr#1#2{\smash{\mathop{\hbox to .5in{\rightarrowfill}}
\limits^{\scriptstyle#1}_{\scriptstyle#2}}}

\def\harrl#1#2{\smash{\mathop{\hbox to .5in{\leftarrowfill}}
\limits^{\scriptstyle#1}_{\scriptstyle#2}}}

\newcommand{\qed}{$\FBox$ \\[1em]}

\newcommand{\FBox}{\rule{2mm}{2.25mm}}
\newcommand{\be}{\begin{equation}}
\newcommand{\ee}{\end{equation}}
\newcommand{\R}{\rm I \mkern -3mu R}

\newtheorem{thm}{Theorem}[section]
\newtheorem{rem}{Remark}[thm]
\newtheorem{lemma}[thm]{Lemma}
\newtheorem{cor}[thm]{Corollary}
\newtheorem{prop}[thm]{Proposition}
\begin{document}
\begin{titlepage}
\begin{center}
{\bf \Large{HIGHER-ORDER LAGRANGIAN FORMALISM ON GRASSMANN MANIFOLDS}}
\end{center}
\vskip 1.0truecm
\centerline{D. R. Grigore
\footnote{e-mail: grigore@theor1.ifa.ro, grigore@roifa.ifa.ro}}
\vskip5mm
\centerline{Dept. of Theor. Phys., Inst. Atomic Phys.}
\centerline{Bucharest-M\u agurele, P. O. Box MG 6, ROM\^ANIA}
\vskip 2cm
\bigskip \nopagebreak
\begin{abstract}
\noindent
The Lagrangian formalism on a arbitrary non-fibrating manifold is considered.
The kinematical description of this generic situation is based on the concept
of (higher-order) Grassmann manifolds which is the factorization of the
regular velocity manifold to the action of the differential group. Here we
introduce in this context the basic concepts of the Lagrangian formalism as
Lagrange, Euler-Lagrange and Helmholtz-Sonin forms.  These objects come in
pairs, namely we have homogeneous objects (defined on the regular velocity
manifold) and non-homogeneous objects (defined on the Grassmann manifold).
We will establish the connection between the homogeneous objects and their
non-homogeneous counterparts. As a result we will conclude that the generic
expressions for a variationally trivial Lagrangian and for a locally
variational differential equation remain the same as in the fibrating case.
\end{abstract}



\end{titlepage}

\section{Introduction}

The Lagrangian formalism is usually based, from the kinematical point of view,
on a fibre bundle with the base manifold and the fibres interpreted as the
``space-time" variables and the field variables respectively. As proved in the
literature, the natural object associated to such a fibre bundle (from the
Lagrangian formalism point of view) is the so-called variational exact
sequence \cite{Kr4} - \cite{Kr6}. This sequence contains as distinguished terms
the Lagrange, Euler-Lagrange and Helmholtz-Sonin forms, which are the main
ingredients of the Lagrangian formalism. Using the exactness property one can
obtain, beside a intrinsic geometrical formulation of the Lagrangian formalism,
the most general expression of a variationally trivial Lagrangian and the
generic form of a locally variational differential equation  \cite{Gr3} -
\cite{Gr4}. 

One of the questions not settled in the literature till recently was if the
formalism above can be extended to the case when the kinematics of the
dynamical system is described by a manifold which is not a fibre bundle. A
physical example of this kind is the Minkowski space describing the 
relativistic particle. For first-order Lagrangian systems a generalization of
the Lagrangian formalism covering this case was proposed (see \cite{GP} and 
\cite{Gr2} for a review) based on the notion of Grassmann manifold and the
so-called Lagrange-Souriau form. For higher-order Lagrangian systems the
construction of the corresponding Grassmann manifold is more subtle and was
performed in \cite{GK}. The idea is to start with the manifold of jets of
immersions in the kinematical manifold of the problem. This manifold is usually
called the velocity manifold and physically corresponds to parametrised
evolutions.  There exists a natural action of the so-called differential group
on this manifold which physically corresponds to changing the parametrisation.
One considers the submanifold of the regular velocities and takes its
factorisation to the differential group. This is exactly the Grassmann manifold
associated to the kinematical manifold of the problem. It was proved that this
is the natural framework to describe a Lagrangian system in this more general
non-fibrating case. The main combinatorial difficulty consists in establishing
a convenient chart system on this factor manifold.  

In this paper we continue the analysis of Lagrangian systems in this framework
by constructing the corresponding Lagrange, Euler-Lagrange and Helmholtz-Sonin
form. The problem to be solved is that the expressions from the fibrating case
are no longer well defined geometrical objects so one must find out proper
substitutes for them. The idea is to construct these kind of objects first on
the velocity manifold as {\it bona fide} geometrical forms and impose some
homogeneity properties. One discovers that these globally defined objects
are inducing {\it locally} defined expressions on the Grassmann manifold which
have convenient transformation properties with respect to a change of charts
and formally coincide with the desired expressions of the usual Lagrangian
formalism. In this way we will be able to define on the Grassmann manifold the
classes (modulo contact forms) of the Lagrange, Euler-Lagrange and
Helmholtz-Sonin forms.

The paper is organised as follows. In Section 2 we remind the basic
construction of a Grassmann manifold following essentially \cite{GK} but
also providing some new results. We do that because we will need many
formul\ae~ for the next sections. For completeness we will also sketch the proof
of the main results. In Section 3 we define the main objects of the Lagrangian
formalism in the non-fibrating case. In Section 4 we give new proofs for the
construction of the Lagrange-Souriau form in the case of first-order Lagrangian
systems and suggests some ways of generalisation of the construction to the
general case of higher-order systems.

\newpage

\section{Grassmann Manifolds}

\subsection{The Manifold of $(r,n)$-Velocities}

Let us consider
$N$, $n\ge 1$
and
$r\ge 0$
integers such that
$n\le N$,
and let
$X$
be a smooth manifold of dimension
$N$ describing the kinematical degrees of freedom of a certain physical
problem.

We will consider 
$U \subset \R^{n}$
a neighbourhood of the point 
$0 \in \R^{n}$,
$x \in X$
and let
$\Gamma_{(0,x)}$
be the set of smooth immersions 
$\gamma: U \rightarrow X$
such that 
$\gamma(0) = x$.
As usual, we consider on
$\Gamma_{(0,x)}$
the relation 
$``\gamma \sim \delta"$
{\it iff} there exists a chart 
$(V,\psi) \quad \psi = (x^{A}), \quad A =1, \dots ,N$
on $X$ such that the functions
$\psi \circ \gamma, \psi \circ \delta: \R^{n} \rightarrow \R^{N}$
have the same partial derivatives up to order $r$ in the point $0$.
One can prove as in \cite{Gr3} that
$\sim$
is a (chart independent) equivalence relation. By an
$(r,n)$-{\it velocity}
at a point
$x \in X$
we mean such an equivalence class of the type
$\Gamma_{(0,x)} /\sim$. 
The equivalence class of
$\gamma$
will be denoted by
$j^{r}_{0}\gamma$.
The set of
$(r,n)$-velocities
at
$x$
is denoted by
$T^{r}_{(0,x)}(\R^{n},Y) \equiv \Gamma_{(0,x)} /\sim$.

Further, we denote
$$
T^{r}_{n}X = \bigcup_{x \in X} T^{r}_{(0,x)}(\R^{n},X), 
$$
and define surjective mappings
$\tau^{r,s}_{n}: T^{r}_{n}X \rightarrow T^{s}_{n}X$, 
where
$0 < s \leq r$,
by
$
\tau^{r,s}_{n}(j^{r}_{0}\gamma) = j^{s}_{0}\gamma
$
and
$\tau^{r,0}_{n}: T^{r}_{n}X \rightarrow X$, 
where
$1 \leq r$,
by
$
\tau^{r,0}_{n}(j^{r}_{0}\gamma) = \gamma(0).
$

In the conditions above let
$(V,\psi), \quad \psi =(x^{A})$,
be a chart on
$X$.
Then we define the couple 
$(V^{r}_{n},\psi^{r}_{n})$
where
$
V^{r}_{n} = (\pi^{r,0}_{n})^{-1}(V), \quad
\psi^{r}_{n} = (x^{A},x^{A}_{j}, \cdots, x^{A}_{j_{1},j_{2},\dots,j_{r}}),
$
where
$
1\leq j_{1} \leq j_{2} \leq \cdots \leq j_{r} \leq n, 
$
and 

\be
x^{A}_{j_{1},\dots,j_{k}}(j^{r}_{0}\gamma) \equiv 
\left. {\partial^{k} \over \partial t^{j_{1}} \dots \partial t^{j_{k}}}
x^{A} \circ \gamma\right|_{0}, 
\quad 0 \le k\le r.
\label{coord-J}
\ee

\begin{rem}
Let us note that the expressions
$x^{A}_{j_{1}, \cdots j_{k}}(j^{r}_{0}\gamma)$
are defined for all values
$j_{1},\dots ,j_{r} \in \{1,\dots ,n\}$
but because of the symmetry property
\be
x^{A}_{j_{P(1)},\dots,j_{P(k)}}(j^{r}_{0}\gamma) =
x^{A}_{j_{1},\dots,j_{k}}(j^{r}_{0}\gamma) \quad (k = 2,...,n)
\ee
for all permutations
$P \in {\cal P}_{k}$
of the numbers
$1, \dots ,k$
we consider only the independent components given by the restrictions
$
1\leq j_{1} \leq j_{2} \leq \cdots \leq j_{r} \leq n.
$
Taking this into account one can use multi-index notations i.e.
$
\psi^{r}_{n} = (x^{A}_{J}), \quad |J| = 0,...,r
$
where by definition
$x^{A}_{\emptyset} \equiv x^{A}.$
The same comment is true for the partial derivatives
${\partial \over \partial x^{A}_{j_{1},\dots,j_{k}}}$.
\end{rem}

Then one can prove that the couple 
$(V^{r}_{n},\psi^{r}_{n})$
is a chart on 
$T^{r}_{n}X$
called the {\it associated chart} of the chart
$(V,\psi)$.
Next, one shows that the set
$T^{r}_{n}X$
has a smooth structure defined by the system of charts
$(V^{r}_{n},\psi^{r}_{n})$;
moreover 
$T^{r}_{n}X$
is a fibre bundle over $X$ with the canonical projection
$\tau^{r,0}$.
The set
$T^{r}_{n}Y$
endowed with the smooth structure defined by the associated charts defined
above is called the {\it manifold of}
$(r,n)$-{\it velocities}
over
$X$.

The equations of the mapping
$\tau^{r,s}_{n}: T^{r}_{n}X \to T^{s}_{n}X$ in terms of the associated
charts are given by 
$
x^{A}_{j_{1},\dots,j_{k}} \circ \tau^{r,s}_{n}(j^{r}_{0}\gamma) =
x^{A}_{j_{1},\dots,j_{k}}(j^{r}_{0}\gamma), 
$
where
$
0\leq k\leq s.
$
These mappings are all submersions. 

\subsection{Formal Derivatives}

Like in \cite{AD}, \cite{An}, \cite{Gr3}, let us consider in the chart
$(V^{r}_{n},\psi^{r}_{n})$
the following differential operators

\be
\partial^{j_{1},\dots,j_{k}}_{A} \equiv {r_{1}! \dots r_{n}! \over k!}
{\partial \over \partial x^{A}_{j_{1},\dots,j_{k}}}, 
\quad j_{1},\dots,j_{k} \in \{1, \dots, n\}
\label{partial}
\ee
where 
$r_{k}$
is the number of times the index $k$ shows up in the sequence
$j_{1}, \dots j_{k}$.

Then one can prove that the following relation is true:

\be
\partial^{i_{1},\dots,i_{k}}_{A} x^{B}_{j_{1},\dots,j_{l}} =
\cases{ \delta^{B}_{A} {\cal S}^{+}_{j_{1},\dots,j_{k}} \delta^{i_{1}}_{j_{1}}
\dots \delta^{i_{k}}_{j_{k}} & if $k = l$ \cr 0 & if $k \not= l$ \cr}.
\label{derivation}
\ee

Here we use the notations from \cite{Gr2}, namely
${\cal S}^{\pm}_{j_{1},\dots,j_{k}}$
is the symmetrization (for the sign $+$) and respectively the 
antisymmetrization (for the sign $-$) projector operator defined by

\be
{\cal S}^{\pm}_{j_{1},\dots,j_{k}} f_{j_{1},\dots,j_{k}} \equiv
{1 \over k!} \sum_{P \in {\cal P}_{k}} \epsilon_{\pm}(P) 
f_{j_{P(1)},\dots,j_{P(k)}}
\label{sa}
\ee
where the sum runs over the permutation group 
${\cal P}_{k}$
of the numbers
$1, \dots, k$
and
$$
\epsilon_{+}(P) \equiv 1, \quad \epsilon_{-}(P) \equiv (-1)^{|P|}, \quad
\forall P \in {\cal P}_{k};
$$
here
$|P|$
is the signature of the permutation $P$.

The differential operators defined by (\ref{partial}) take care of overcounting
the indices. More precisely, for any smooth function on 
$V^{r}$,
the following formula is true:

\be
df = \sum_{k=0}^{r} (\partial^{j_{1},\dots,j_{k}}_{A} f) 
dx^{A}_{j_{1},\dots,j_{k}} = 
\sum_{|I| \leq r} (\partial^{I}_{A} f) dx^{A}_{I}
\label{df}
\ee
where we have also used the convenient multi-index notation.

We define now in the chart
$(V^{r}_{n},\psi^{r}_{n})$
the {\it formal derivatives} by the expressions

\be
d^{r}_{i} \equiv \sum_{k=0}^{r-1} x^{A}_{i,j_{1},\dots,j_{k}} 
\partial^{j_{1},\dots,j_{k}}_{A} = 
\sum_{|J| \leq r-1} x^{A}_{iJ} \partial^{J}_{A}.
\label{formal}
\ee

The last expression uses the multi-index notation; if $I$ and $J$ are two such
multi-indices we mean by $IJ$ the juxtaposition of the two sets $I, J.$

We note that the preceding formula does not define a vector field on 
$T^{r}_{n}Y$.
When no danger of confusion exists we simplify the notation putting simply
$d_{i} = d^{r}_{i}.$
One can easily verify that the following formul\ae\/ follow directly from the
definition: 

\be
d_{i} x^{A}_{j_{1},\dots,j_{k}} = \cases{x^{A}_{i,j_{1},\dots,j_{k}} & 
if $k \leq r-1$ \cr 0 & if $k = r$ \cr},
\label{d-x}
\ee

\be
\left[ \partial^{j_{1},\dots,j_{k}}_{A}, d_{i} \right] =
{\cal S}^{+}_{j_{1},\dots,j_{k}} \delta^{j_{1}}_{i} 
\partial^{j_{2},\dots,j_{k}}_{A}, \quad k = 0,...,r
\ee
and
\be
\left[ d_{i}, d_{j} \right] = 0.
\label{commutator}
\ee

The formal derivatives can be used to conveniently express the change of charts
on the velocity manifold induced by a change of charts on $X$. Let
$(V,\psi)$
and
$(\bar{V},\bar{\psi})$
two charts on $X$ such that 
$V \cap \bar{V} \not= \emptyset$
and let
$(V^{r},\psi^{r})$
and
$(\bar{V}^{r},\bar{\psi}^{r})$
the corresponding attached charts from
$T^{r}_{n}X.$
The change of charts on $X$ is
$F: \R^{N} \rightarrow \R^{N}$
given by:
$F \equiv \bar{\psi} \circ \psi^{-1}.$
It is convenient to denote by 
$F^{A}: \R^{N} \rightarrow \R$
the components of $F$ given by
$F^{A} \equiv \bar{x}^{A} \circ \psi^{-1}$.
We now consider the change of charts on 
$T^{r}_{n}X$
given by
$F^{r} \equiv \bar{\psi}^{r} \circ (\psi^{r})^{-1}.$
One notes that
$V^{r} \cap \bar{V}^{r} \not= \emptyset$;
we need the explicit formul\ae\/ for the components of 
$F^{r}$, namely for the functions
$$
F^{A}_{j_{1}, \dots, j_{k}} \equiv \bar{x}^{A}_{j_{1}, \dots, j_{k}} \circ
(\psi^{r})^{-1}, \quad j_{1} \leq j_{2} \dots \leq j_{k}, \quad k = 1,...,r
$$
defined on the overlap:
$V^{r} \cap \overline V^{r}$.
First one notes the following relation:

\be
\overline d_{i} = d_{i}.
\label{dd}
\ee

Indeed, one defines for any immersion
$\gamma \in \Gamma_{(0,x)}$
the map
$j^{r}\gamma$
from
$\R^{n}$
into
$T^{r}_{n}X$
given by

\be
x^{A}_{j_{1},\dots,j_{k}} \circ j^{r}\gamma(t) \equiv 
{\partial^{k} x^{A} \circ \gamma \over 
\partial t^{j_{1}} \dots \partial t^{j_{k}}}(t)
\quad 0 \le k\le r
\label{extension}
\ee
and easily discovers that

\be
(j^{r}\gamma)_{*0}{\partial \over \partial t^{i}} = d_{i} = \bar{d}_{i}.
\ee

Using (\ref{dd}) one easily finds out that the functions
$F^{A}_{j_{1}, \dots, j_{k}}$
are given recurringly by the following relation:

\be
F^{A}_{jI} = d_{j} F^{A}_{I} \quad |I| \leq r-1;
\label{F_I}
\ee
(compare with (\ref{d-x}).)

This relation can be ``solved" explicitly according to

\begin{lemma}
The following formula holds

\be
F_{I}^{A} = \sum_{p=1}^{|I|} 
{\sum\limits_{(I_{1},\ldots,I_{p})}} 
x_{I_{1}}^{B_{1}} \cdots x_{I_{p}}^{B_{p}}
(\partial_{B_{1}} \cdots \partial_{B_{p}} F^{A}), \quad 1 \leq |I| \leq r
\label{change-chart}
\ee
where the second sum denotes summation over all partitions
${\cal P}(I)$
of the set
$I$
and two partitions are considered identical if they differ only by a
permutation of the subsets. 
\end{lemma}

{\bf Proof:}
We sketch the proof because the argument will be used repeatedly in this paper. 
It is natural to use complete induction on 
$|I|$. 
For 
$I = \{ j \}$
the formula from the statement coincides with (\ref{F_I}) for 
$I = \emptyset$.
We suppose the formula true for any multi-index $I$ with
$|I| = s < r$
and prove it for the multi-index
$jI$. If we use (\ref{F_I}) we get:

\begin{eqnarray}
F^{A}_{jI} = \sum_{p=1}^{|I|} {\sum\limits_{(I_{1},\ldots,I_{p})}} 
\left[ \sum_{l=1}^{p} x_{I_{1}}^{B_{1}} \cdots (d_{j} x_{I_{l}}^{B_{l}}) \cdots
x_{I_{p}}^{B_{p}}
(\partial_{B_{1}} \cdots \partial_{B_{p}} F^{A}) +
x_{I_{1}}^{B_{1}} \cdots x_{I_{p}}^{B_{p}}
d_{j} (\partial_{B_{1}} \cdots \partial_{B_{p}} F^{A}) \right] \nonumber \\ 
= \sum_{p=1}^{|I|} {\sum\limits_{(I_{1},\ldots,I_{p})}} 
\left[ \sum_{l=1}^{p} x_{I_{1}}^{B_{1}} \cdots x_{jI_{l}}^{B_{l}} \cdots
x_{I_{p}}^{B_{p}}
(\partial_{B_{1}} \cdots \partial_{B_{p}} F^{A}) +
x_{I_{1}}^{B_{1}} \cdots x_{I_{p}}^{B_{p}}
x^{B_{p+1}}_{j} (\partial_{B_{1}} \cdots \partial_{B_{p+1}} F^{A}) \right].
\nonumber
\end{eqnarray}

We now note that the partitions 
${\cal P}(jI)$
of the set $jI$
can be obtained in two distinct ways:
\begin{itemize}
\item
by taking a partition 
$(I_{1},\ldots,I_{p}) \in {\cal P}(I)$
and adjoining the index $j$ to
$I_{1}, I_{2},\dots,I_{p}$;
\item
by taking a partition
$(I_{1},\ldots,I_{p}) \in {\cal P}(I)$
and constructing the associated partition \\
$(I_{1},\ldots,I_{p},j) \in {\cal P}(jI).$
\end{itemize}

We get the two types of contributions in the formula above and this finishes
the proof.
$\qed$

\begin{rem}
The combinatorial argument above will be called {\bf the partition argument}.
\end{rem}

\begin{rem}
From the formula derived above it immediately follows that we have:

\be
\partial^{J}_{B} F^{A}_{I} = 0,  \quad 0 \leq |I| < |J| \leq r
\ee
i.e. the functions
$F^{A}_{I}$
depend only of the variables
$x^{B}_{J}$
with the restrictions specified above.
\end{rem}

\subsection{The Differential Group}

By definition the {\it differential group of} order $r$ is the set

\be
L^{r}_{n} \equiv \{ j^{r}_{0} \alpha \in J^{r}_{0,0}(\R^{n},\R^{n}) \vert
\alpha \in Diff(\R^{n}) \}
\ee
i.e. the group of invertible
$r$-jets
with source and target at
$0\in \R^{n}$.
The group multiplication in
$L^{r}_{n}$
is defined by the jet composition
$
L^{r}_{n}\times L^{r}_{n} \ni (j^{r}_{0}\alpha, j^{r}_{0}\beta) 
\mapsto  j^{r}_{0}(\alpha\circ \beta) \in L^{r}_{n}.
$

The {\it canonical} (global) {\it coordinates} on 
$L^{r}_{n}$
are defined by

\be
a^{i}_{j_{1},\dots,j_{k}}(j^{r}_{0}\alpha) = \left.
{\partial^{k} \alpha^{i} \over 
\partial t^{j_{1}} \dots \partial t^{j_{k}}}\right|_{0}, \quad
j_{1} \leq j_{2} \leq \dots \leq j_{k},  \quad k = 0,...,r
\label{coord-L}
\ee
where
$\alpha^{i}$
are the components of a representative
$\alpha$
of
$j^{r}_{0}\alpha$.

We denote
$$
a \equiv (a^{i}_{j},a^{i}_{j_{1},j_{2}}, \dots ,a^{i}_{j_{1},\dots,j_{k}}) = 
(a^{i}_{J})_{|J| \leq r}
$$
and notice that one has

\be
det(a^{i}_{j}) \not= 0.
\label{group}
\ee 

To obtain the composition law for the differential group we need a
combinatorial result following easily by induction with the partition argument:

\begin{lemma}
Let 
$U, V \in \R^{n}$
be open sets,
$\alpha: U \rightarrow V$
and
$f: V \rightarrow \R$
smooth functions. Then the following formula is true:

\be
\partial_{I} (f \circ \alpha) = \sum_{p=1}^{|I|} 
{\sum\limits_{(I_{1},\ldots,I_{p})}} 
(\partial_{I_{1}} \alpha^{i_{1}}) \dots (\partial_{I_{p}} \alpha^{i_{p}})
(\partial_{i_{1},\dots ,i_{p}} f) \circ \alpha
\ee
where we have denoted for any multi-index
$I = \{ i_{1}, \dots , i_{s}\}$
$$
\partial_{I} f \equiv {\partial^{s} f \over \partial t^{i_{1}} \dots
\partial t^{i_{s}}}.
$$
\label{der-prod}
\end{lemma}

We now have:

\begin{lemma}
The group multiplication in
$L^{r}_{n}$
is expressed in the canonical coordinates by the equations 

\be
(a \cdot b)^{k}_{I} = \sum^{|I|}_{p=1} {\sum\limits_{(I_{1},\dots,I_{p})}}
b^{j_{1}}_{I_{1}} \dots b^{j_{p}}_{I_{p}} a^{k}_{j_{1}, \dots j_{p}}, 
\quad |I| = 1, \dots r. 
\label{composition}
\ee

The group
$L^{r}_{n}$
is a Lie group.
\end{lemma}

{\bf Proof:}
(i) We start from the defining formula:
$$
(a\cdot b)^{k}_{j_{1},\dots,j_{l}} = \left. 
{\partial^{l} \alpha^{k} \circ \beta \over 
\partial t^{j_{1}} \dots \partial t^{j_{l}}}\right|_0
$$
and apply the lemma above. One obtains the composition formula.

(ii) It is clear that the composition formula (\ref{composition}) is a smooth
function. The identity is evidently:
$$
e \equiv (\delta^{i}_{j},0, \dots ,0)
$$
and it remains to prove that the map
$a \rightarrow a^{-1}$
is smooth. Indeed one immediately proves by induction that
$$
(a^{-1})^{k}_{I} = \left( det(a^{i}_{j}) \right)^{-|I|} \times
P^{k}_{I}(a)
$$
where 
$P^{k}_{I}$
is a polynomial in the variables
$a^{i}_{I}, \quad |I| = 0, \dots ,r.$
$\qed$

The manifolds of
$(r,n)$-velocities
$T^{r}_{n}Y$
admits a (natural) smooth right action of the differential group
$L^{r}_{n}$,
defined by the jet composition

\be
(x\cdot a)^{A}_{I} \equiv x^{A}_{I}(j_{0}^{r}(\gamma \circ \alpha)) 
\label{x-a}
\ee
where the connection between 
$x^{A}_{I}$
and
$\gamma$
is given by (\ref{coord-J}) and the connection between 
$a^{i}_{I}$
and
$\alpha$
is given by (\ref{coord-L}). 

We determine the chart expression of this action. 

\begin{prop}
The group action (\ref{x-a})is expressed by the equations 

\be
(x\cdot a)^{A} = x^{A},\quad 
(x\cdot a)_{I}^{A} = 
\sum_{p=1}^{|I|} {\sum\limits_{(I_{1},\dots,I_{p}) \in {\cal P(I)}}} 
a_{I_{1}}^{j_{1}} \dots a_{I_{p}}^{j_{p}} x_{j_{1},\dots,j_{p}}^{A},
\quad |I| \geq 1
\label{action}
\ee
and it is smooth.
\end{prop}

{\bf Proof:}

One applies the definitions (\ref{coord-J}) and (\ref{coord-L}) together with
the lemma \ref{der-prod}. The smoothness is obvious from the explicit action
formula given above. 
$\qed$

The group
$L^{r}_{n}$
has a natural smooth left action on the set of smooth real functions defined on
$T^{r}_{n}X$ ,
namely for any such function $f$ we have:

\be
(a\cdot f)(x) \equiv f(x\cdot a).
\label{act}
\ee

\subsection{Higher Order Regular Velocities}

We say that a
$(r,n)$-velocity
$j_{0}^{r}\gamma \in T_{n}^{r}X$
is {\it regular}, if 
$\gamma$
(or any other representative) is an immersion. If
$(V,\psi)$,
$\psi =(x^{A})$,
is a chart, and the target
$\gamma(0)$
of an element
$j_{0}^{r}\gamma \in T_{n}^{r}X$
belongs to
$V$,
then
$j_{0}^{r}\gamma$
is regular {\it iff} there exists a subsequence 
$ {\bf I} \equiv (i_{1},\dots,i_{n})$ 
of the sequence
$(1,2,\dots,n,n+1,\dots,n+m)$
such that

\be
{\rm det} (x^{i_{k}}_{j}) \not= 0;
\label{regular}
\ee
(here 
$x^{i_{k}}_{j}$
is a 
$n \times n$
real matrix.)

The associated charts have the form
$$
(V^{{\bf I},r},\psi^{{\bf I},r}), \quad
\psi^{{\bf I},r} = (x^{k}_{I},x^{\sigma}_{I}), \quad k = 1,\dots,n, 
\quad \sigma = 1,\dots m \equiv N - n, \quad |I| \leq r
$$
where
$$
x^{k}_{I} \equiv x^{i_{k}}_{I}, \quad k = 1,\dots n
$$ 
and
$\sigma \in \{ 1,\dots, N\} - \{i_{1},\dots,i_{n}\}.$
The set of regular
$(r,n)$-velocities
is an open,
$L_{n}^{r}$-invariant
subset of
$T_{n}^{r}X$,
which is called the {\it manifold of regular} $(r,n)$-{\it velocities},
and is denoted by
${\rm Imm} T_{n}^{r}X$.

We want to find out a complete
system of
$L_{n}^{r}$-invariants
(in the sense of Weyl) of the action (\ref{action}) on
${\rm Imm} T_{n}^{r}X$.

We will consider, for simplicity a chart for which one has
$\{i_{1},\dots,i_{n}\} = \{ 1,\dots, n\}$
and we will denote
$$
x^{\sigma}_{I} \equiv x^{n+\sigma}_{I}, \quad 
\sigma = 1,\dots m, \quad |I| \leq r.
$$

We begin with the following result:

\begin{prop}
Let
$(x^{\sigma}_{I},x^{i}_{I})$
be the coordinates of a point in
${\rm Imm} T_{n}^{r}X$.
Then

\be
{\bf x} \equiv (x^{i}_{I})_{1 \leq |I| \leq r}
\label{x}
\ee
is a element from
$L_{n}^{r}$ . We denote its inverse by

\be
{\bf z} \equiv (z^{i}_{I})_{1 \leq |I| \leq r}.
\ee

Then
$z^{i}_{j}$
is the inverse of the matrix
$x^{l}_{p}$:

\be
z^{i}_{j} x^{j}_{p} = \delta^{i}_{p}
\ee
and the functions
$z^{i}_{j_{1},\dots,j_{k}}, k = 2,\dots r$
can be determined recurringly from the equations:

\be
z_{j_{1},\dots,j_{k}}^{i} = z_{j_{1}}^{p} d_{p}
z_{j_{2},\dots,j_{k}}^{i},\quad k = 2,\dots, r
\label{z}
\ee
\end{prop}

{\bf Proof:}
For the first assertion one uses (\ref{group}) and (\ref{regular}). For the
relation (\ref{z}) one starts from the definition
$
{\bf z} \cdot {\bf x} = e
$
or, in detail

$$
\sum^{|I|}_{k=1} {\sum\limits_{(I_{1},\dots,I_{k})}}
x^{j_{1}}_{I_{1}} \dots x^{j_{k}}_{I_{k}} z^{i}_{j_{1},\dots,j_{k}} =
\cases{ \delta^{i}_{I} & {\rm for} $|I| = 1$ \cr 0 & {\rm for} 
$|I| = 2, \dots r$ \cr}.
$$

One performs two distinct operations on this relation: (a) we apply the 
operator
$d_{p}$;
(b) we make
$I \mapsto Ip$.
Next one subtracts the two results and uses the partition argument;  the
following relation follows:

$$
\sum^{|I|}_{k=1} {\sum\limits_{(I_{1},\dots,I_{k})}}
x^{j_{1}}_{I_{1}} \dots x^{j_{k}}_{I_{k}} (d_{p} z^{i}_{j_{1},\dots,j_{k}} -
x^{j_{0}}_{p} z^{i}_{j_{0},\dots,j_{k}} ) = 0.
$$
 
From this relation, we obtain, by induction the formula from the statement.
$\qed$

The formula (\ref{z}) suggests the following result:

\begin{prop}

Let
$(V,\psi), \psi = (x^{A})$,
be a chart on
$X$
and let
$(V_{n}^{r},\psi_{n}^{r})$
be the associated chart on
${\rm Imm} T_{n}^{r}X$.
We define recurringly on this chart the following functions

\be
y^{\sigma} \equiv x^{\sigma}, \quad
y_{i_{1},\dots, i_{k}}^{\sigma} = z_{i_{1}}^{j} d_{j}
y_{i_{2},\dots, i_{k}}^{\sigma},\quad k = 1,\dots, r;
\label{invariants}
\ee
(here 
$z^{j}_{i}$
are the first entries of the element 
${\bf z} \in L^{r}_{n}$.)

Then the functions
$y^{\sigma}_{i_{1},\dots,i_{k}}$
so defined depend smoothly only on
$x^{A}_{J}, \quad |J| \leq k$
and are completely symmetric in all indices
$i_{1},\dots, i_{k}, \quad k = 1,\dots r.$
\end{prop}

{\bf Proof:}

The first assertion follows immediately by induction. Next, one derives
directly from the formula (\ref{invariants}) that

$$
y^{\sigma}_{i_{1},\dots,i_{k}} = z^{j_{1}}_{i_{1}} z^{j_{2}}_{i_{2}}
\left( d_{j_{1}} d_{j_{2}} y^{\sigma}_{i_{3},\dots,i_{k}} -
x^{p}_{i_{1},i_{2}} z^{j}_{p} d_{j} y^{\sigma}_{i_{3},\dots,i_{k}}\right),
\quad k = 2,\dots r.
$$

In particular we see that for $k = 2$ the symmetry property is true. One can
proceed now by induction. If
$y^{\sigma}_{i_{1},\dots,i_{k-1}}$
is completely symmetric then the formula above shows that we have the symmetry
property in the indices
$i_{1}$ 
and 
$i_{2}$; 
moreover the recurrence relation (\ref{invariants}) shows that we have the
symmetry property in the indices
$i_{2},\dots, i_{k}$.
So we obtain the desired property in all indices. 
$\qed$

As a result of the symmetry property just proved we can use the convenient
multi-index notation
$y^{\sigma}_{I}, \quad |I| \leq r.$
Now we have an explicit formula for these functions

\begin{prop}
The functions
$y^{\sigma}_{I}, \quad 1 \leq |I| \leq r$
are {\bf uniquely} determined by the recurrence relations:

\be
x_{I}^{\sigma} = \sum_{p=1}^{|I|} \sum\limits_{(I_{1},\dots,I_{p})}
x_{I_{1}}^{j_{1}} \dots x_{I_{p}}^{j_{p}} y_{j_{1},\dots,j_{p}}^{\sigma}.
\label{yx}
\ee

Using the notation 
${\bf x} \in L^{r}_{n}$
one can compactly write the relation above as

\be
x^{\sigma}_{I} = (y\cdot {\bf x})^{\sigma}_{I}, \quad 1 \leq |I| \leq r.
\label{yx-compact}
\ee
\end{prop}

{\bf Proof:}
Goes by induction on 
$|I|$.
The formula above is obvious for
$I = \{ j\}$.
If it is valid for 
$|I| < r$
we apply to the relation above the operator
$d_{j}$
and use (\ref{d-x}) and the partition argument. One obtains in this way the
formula from the statement for
$Ij$.
The unicity also follows by induction. The last assertion is a consequence of
the first formula and of the expression of the group action (\ref{action}).
$\qed$

Let us note that one can ``invert" the formul\ae\/ from the statement. Indeed,
(\ref{yx-compact}) is equivalent to

\be
y^{\sigma}_{I} = (x\cdot {\bf z})^{\sigma}_{I}, \quad 1 \leq |I| \leq r
\label{xy-compact}
\ee
or, explicitly:

\be
y_{I}^{\sigma} = \sum_{p=1}^{|I|} \sum\limits_{(I_{1},\dots,I_{p})}
z_{I_{1}}^{j_{1}} \dots z_{I_{p}}^{j_{p}} x_{j_{1},\dots,j_{p}}^{\sigma}.
\label{xy}
\ee

\begin{cor}

One can use on
$V^{r}$
the new coordinates
$(y^{\sigma}_{I},x^{i}_{I}), \quad |I| \leq r$.
\label{new}
\end{cor}

{\bf Proof:} From the relations (\ref{yx}) and (\ref{xy}). 
$\qed$

Now we have the following result

\begin{prop}

The functions
$y^{\sigma}_{I}, \quad |I| \leq r$
are 
$L^{r}_{n}$-invariants with respect to the natural action (\ref{act}).
\end{prop}

{\bf Proof:}
Let 
$a \in L^{r}_{n}$
be arbitrary. We start from (\ref{yx-compact}) 
and use the associativity of the group composition law of
$L^{r}_{n}$;
we get:

$$
(x\cdot a)^{\sigma}_{I} =  ((y\cdot {\bf x})\cdot a)^{\sigma}_{I} = 
(y \cdot ({\bf x}\cdot a))^{\sigma}_{I}.
$$

On the other hand if we make in (\ref{yx-compact}) the substitution 
$x \mapsto x\cdot a$
we get: 

$$
(x\cdot a)_{I}^{\sigma} = ((a\cdot y)\cdot (x\cdot a))_{I}^{\sigma}.
$$
(here 
$a\cdot y$
denotes the action of the differential group on the functions $y$ according to
(\ref{act}).)

By comparing the two formul\ae\/ and using of the unicity statement from the
preceding proposition we get the desired result.
$\qed$

Moreover, we can prove:

\begin{thm}

The functions
$y^{\sigma}_{I}, \quad |I| \leq r$
are a complete system of invariants in the sense of Weyl.
\end{thm}

{\bf Proof:}
The fact that the functions
$y^{\sigma}_{I}, \quad |I| \leq r$
are functionally independent follows by {\it reductio ad absurdum}. If they
would be functionally dependent, then from (\ref{xy}) it would follow that the
expressions 
$x^{\sigma}_{I}, \quad |I| \leq r$ 
also are functionally dependent.

We still must show that there are no other invariants beside
$y^{\sigma}_{I}, \quad |I| \leq r$.
We proceed as follows. From corollary \ref{new} it follows that one can use on
$V^{r}$
the coordinates
$(y^{\sigma}_{I},x^{i}_{I}),  \quad |I| \leq r$.
In these coordinates the action of the group
$L^{r}_{n}$ 
is:

\begin{eqnarray}
(y\cdot a)^{\sigma}_{I} = y^{\sigma}_{I}, \quad |I| \leq r, \nonumber \\
(x \cdot a)^{i} = x^{i}, \quad
(x \cdot a)^{i}_{I} = \sum_{p=1}^{|I|} \sum\limits_{(I_{1},\dots,I_{p})}
a_{I_{1}}^{j_{1}} \dots a_{I_{p}}^{j_{p}} x_{j_{1},\dots,j_{p}}^{i},
\quad 1 \leq |I| \leq r.
\label{new-action}
\end{eqnarray}

One can prove now by induction that this action is transitive. This shows that
the system of invariants from the statement is complete.
$\qed$

\subsection{Higher Order Grassmann Bundles}

Needless to say, the whole formalism presented above can be implemented in an
arbitrary chart system 
$(V^{{\bf I},r},\psi^{{\bf I},r})$
on
${\rm Imm}T^{r}_{n}X$
(see the beginning of the preceding subsection). In this context we finally 
have the central result:

\begin{thm}

The set
$
P^{r}_{n}X \equiv {\rm Imm} T_{n}^{r}X / L^{r}_{n}
$
has a unique differential manifold structure such that the canonical projection
$\rho^{r}_{n}$ 
is a submersion. The group action (\ref{action}) defines on
${\rm Imm} T_{n}^{r}X$
the structure of a right principal
$L_{n}^{r}$-bundle.
 
A chart system on
$
P^{r}_{n}X
$
adapted to this fibre bundle structure is formed from couples
$(W^{{\bf I},r},\Phi^{{\bf I},r})$
where:

\be
W^{{\bf I},r} = \left\{ j_{0}^{r}\gamma \in V^{r}| 
{\rm det}(x_{j}^{i_{k}}(j_{0}^{r}\gamma)) \not= 0 \right\}
\ee
and

\be
\Phi^{{\bf I},r} = (x^{i}_{I}, y^{\sigma}_{I}), \quad |I| \leq r.
\ee

In this case the local expression of the canonical projection is
$$
\rho^{r}_{n}(x^{i}_{I}, y^{\sigma}_{I}) = (x^{i},y^{\sigma}_{I}).
$$
\end{thm}

{\bf Proof:} 
We define on
${\rm Imm} T_{n}^{r}X \times {\rm Imm} T_{n}^{r}X$
the equivalence relation
$$
x \sim \bar{x} \quad {\it iff} \quad \exists a \in L^{r}_{n} \quad s.t. \quad
\bar{x} = x\cdot a.
$$

To prove the first assertion from the statement is sufficient (according to
\cite{D}, par. 16.19.3) to prove that the graph of 
$\sim$
is a closed submanifold of the product manifold. We will look for a convenient
system of coordinates on
${\rm Imm} T_{n}^{r}X$.

The first step is to take
$x$
and
$\bar{x}$
such that
$
x \sim \bar{x}
$
and
to solve the system of equations

$$
\bar{x} = x\cdot a_{x,\bar{x}}
$$
for the unknown functions
$
a_{x,\bar{x}} \in L^{r}_{n}.
$

One easily gets

$$
(a_{x,\bar{x}})^{i}_{j} = z^{i}_{p} \bar{x}^{p}_{j}
$$
and then shows by induction that
$
(a_{x,\bar{x}})^{i}_{I}
$
are uniquely determined smooth functions of
$x^{i}_{J}, \bar{x}^{i}_{J}, \quad |J| \leq |I| \leq r.$

We now define the (local) smooth functions on
${\rm Imm} T_{n}^{r}X \times {\rm Imm} T_{n}^{r}X$:

$$
\Phi^{\sigma}_{I}(x,\bar{x}) \equiv \bar{x}^{\sigma}_{I} - 
(x\cdot a_{x,\bar{x}})^{\sigma}_{I}, \quad |I| \leq r.
$$

It is clear that on can take on
${\rm Imm} T_{n}^{r}X \times {\rm Imm} T_{n}^{r}X$
the (local) coordinates
$(x^{A}_{I}, \Phi^{\sigma}_{I},\bar{x}^{i}_{I}), \quad |I| \leq r$;
it follows that the graph of
$\sim$
is given by

$$
\Phi^{\sigma}_{I} = 0, \quad |I| \leq r
$$
i.e. it is a closed submanifold.

To prove the fibre bundle structure it is sufficient to show (see also 
\cite{D}) that the action of 
$L^{r}_{n}$
is free i.e.
$$
x\cdot a = x \Longrightarrow a = e.
$$
This fact follows elementary by induction.

Finally, we have remarked before (see the preceding theorem) that one can take 
on
${\rm Imm} T_{n}^{r}X$
the coordinates
$(y^{\sigma}_{I},x^{i}_{I}), \quad |I| \leq r$
with the action given by (\ref{new-action}). The last assertion about the
expression of the canonical projection follows.
$\qed$

A point of
$P^{r}_{n}X$
containing a regular
$(r,n)$-velocity
$j^{r}_{0}\gamma$
is called an
$(r,n)$-{\it contact element}, or an
$r$-{\it contact element} of an
$n$-dimensional
submanifold of
$X$,
and is denoted by
$[j^{r}_{0}\gamma].$
As in the case of
$r$-jets, the point
$0 \in \R^{n}$
(resp.
$\gamma(0) \in X$)
is called the {\it source} (resp. the {\it target}) of 
$[j^{r}_{0}\gamma].$
The manifold
$P^{r}_{n}$
is called the
$(r,n)$-{\it Grassmannian bundle},
or simply a {\it higher order Grassmannian bundle} over
$X$.

Besides the quotient projection
$\rho^{r}_{n}: {\rm Imm} T^{r}_{n}X \to P^{r}_{n}$ 
we have for every 
$1 \leq s \leq r$,
the {\it canonical projection} of
$P^{r}_{n}X$
onto
$P^{s}_{n}X$
defined by
$\rho^{r,s}_{n}([j^{r}_{0}\gamma]) = [j^{s}_{0}\gamma]$ 
and the {\it canonical projection} of
$P^{r}_{n}X$
onto
$X$
defined by
$\rho^{r}_{n}([j^{r}_{0}\gamma]) = \gamma(0)$. 

\begin{rem}
When the manifold $X$ is fibred over a manifold $M$ of dimension $n$ one
can also construct the jet extension
$J^{r}X$
(see \cite{Kr4}, \cite{Gr3}). One can establish a canonical isomorphism between
$P^{s}_{n}X$
and
$J^{r}X$
as follows: let 
$x \in M, x = \gamma(m)$,
$\gamma \in \Gamma_{(m, x)}$;
and 
$\phi: \R^{n} \rightarrow M$
a (local) diffeomorphism such that
$\phi(0) = m.$
We can define
$\tilde{\gamma} \in \Gamma_{(0, x)}$
by the formula
$\tilde{\gamma} \equiv \gamma \circ \phi$.
One notices that
$\tilde{\gamma}_{1} \sim \tilde{\gamma}_{2}$
{\it iff} there exists
$\alpha \in Diff(\R^{n})$
such that
$\gamma_{1} = \gamma_{2} \circ \alpha.$
This means that the map
$j^{r}_{n}\gamma \mapsto j^{r}\tilde{\gamma}$
can be factorized to a map from
$P^{r}_{n}X  \rightarrow J^{r}X$
which is proved to be an isomorphism.
\end{rem}

We also note the following result:

\begin{prop}

The following formula is true

\be
\Phi^{\sigma}_{I} = \sum_{p=1}^{|I|} \sum\limits_{(I_{1},\dots,I_{p})}
\bar{x}_{I_{1}}^{j_{1}} \dots \bar{x}_{I_{p}}^{j_{p}} 
(\bar{y}_{j_{1},\dots,j_{p}}^{\sigma} - y_{j_{1},\dots,j_{p}}^{\sigma}),
\quad 1 \leq |I| \leq r.
\label{phi}
\ee
or, in compact notations

\be
\Phi^{\sigma}_{I} = ((\bar{y} - y)\cdot \bar{\bf x})^{\sigma}_{I}.
\label{phi-compact}
\ee

In particular, the equation
$$
\Phi^{\sigma}_{I} = 0, \quad 1 \leq |I| \leq r
$$
is equivalent to
$$
\bar{y}^{\sigma}_{I} = y^{\sigma}_{I}, \quad |I| \leq r.
$$
\end{prop}

{\bf Proof:}
The proof relies heavily on induction. Firstly, we define on
${\rm Imm} T_{n}^{r}X \times {\rm Imm} T_{n}^{r}X$
the expressions

$$
V_{i} \equiv \bar{d}_{i} - (a_{x,\bar{x}})^{i}_{j} d_{j}
$$
(where
$a_{x,\bar{x}}$
have been defined previously) and we prove by induction the following formula:

$$
V_{i} (a_{x,\bar{x}})^{j}_{I} = (a_{x,\bar{x}})^{j}_{iI}, 
\quad 1 \leq |I| \leq r-1.
$$

Next, one uses this formula to prove by direct computation that

$$
V_{i} \Phi^{\sigma}_{I} = \Phi^{\sigma}_{iI}.
$$

Finally, one uses the preceding formula to prove by induction the formula
(\ref{phi}) from the statement.
$\qed$

\begin{rem}
The dimension of the factor manifold
$P_{n}^{r}X$
is

$$
{\rm dim} P_{n}^{r}X = m {n+r \choose n} + n. 
$$
\end{rem}

Now we try to define on
$P_{n}^{r}X$
the analogue of the total differential operators.

\begin{prop}

Let us consider on the regular velocities manifold 
${\rm Imm} T_{n}^{r}X$
the coordinates
$(y^{\sigma}_{I},x^{i}_{I}), \quad |I| \leq r$
and define the operators

\be
\widetilde{\Delta}^{j_{1},\dots j_{k}}_{\sigma} \equiv 
{r_{1}! \dots r_{n}! \over k!}
{\partial \over \partial y^{\sigma}_{j_{1},\dots,j_{k}}}
\label{tilde-partial}
\ee
(where we use the same conventions as in (\ref{partial})).

We also define, by analogy to (\ref{formal})

\be
\tilde{d}^{r}_{i} \equiv {\partial \over \partial x^{i}} + 
\sum_{k=0}^{r-1} y^{\sigma}_{i,j_{1},\dots,j_{k}} 
\widetilde{\Delta}^{j_{1},\dots,j_{k}}_{\sigma} = 
\sum_{|J| \leq r-1} y^{\sigma}_{iJ} \widetilde{\Delta}^{J}_{\sigma}.
\ee

Then the following formula is true

\be
d_{i} = \sum_{p=1}^{r-1} x^{l}_{i,j_{1},\dots,j_{k}} 
\partial^{j_{1},\dots,j_{k}}_{l} + x^{p}_{i} \tilde{d}_{p} 
\label{tilde-d}
\ee
where, as usual,
$\tilde{d}_{i} = \tilde{d}^{r}_{i}$
when no danger of confusion arises.
\end{prop}

{\bf Proof:}
By direct computation.
$\qed$

Now we define on
$P^{r}_{n}X$,
in the chart
$\rho^{r}_{n}(W^{{\bf I},r})$
some operators which are the analogues of (\ref{tilde-partial}) and 
(\ref{formal}), namely

\be
\Delta^{j_{1},\dots j_{k}}_{\sigma} \equiv {r_{1}! \dots r_{n}! \over k!}
{\partial \over \partial y^{\sigma}_{j_{1},\dots,j_{k}}}
\label{delta}
\ee

and

\be
D_{i} \equiv {\partial \over \partial x^{i}} + 
\sum_{k=0}^{r-1} y^{\sigma}_{i,j_{1},\dots,j_{k}} 
\Delta^{j_{1},\dots,j_{k}}_{\sigma} = 
{\partial \over \partial x^{i}} + 
\sum_{|J| \leq r-1} y^{\sigma}_{iJ} \Delta^{J}_{\sigma}.
\label{D}
\ee

These operators are also called {\it total derivatives}. A formula similar to
(\ref{derivation}) is valid and moreover, the preceding proposition has the 
following consequence:

\begin{prop}
The following formula is true:

\be
(\rho^{r}_{n})_{*} (z^{j}_{i} d_{j}) = D_{i}.
\label{Dd}
\ee

In particular, we have for any smooth function $f$ on 
$\rho^{r}_{n}(W^{r})$
the following formula:

\be
d_{i}(f \circ \rho^{r}_{n}) = x^{j}_{i} (D_{j}f) \circ \rho^{r}_{n}.
\label{dD}
\ee

Therefore, if
$(V,\psi)$
and
$(\bar{V},\bar{\psi})$
are two charts on $X$ such that
$V \cap \bar{V} \not= \emptyset$
and 
$D_{i}, \bar{D}_{i}, \quad i = 1,\dots n$
are the corresponding operators defined on
$\rho^{r}_{n}(V^{r})$
and respectively on
$\rho^{r}_{n}(\bar{V^{r}})$,
then we have on
$\rho^{r}_{n}(V^{r} \cap \bar{V^{r}})$:

\be
{\rm Span} (D_{1},\dots,D_{n}) = {\rm Span} (\bar{D}_{1},\dots,\bar{D}_{n}).
\ee
\end{prop}

{\bf Proof:}
The first formula follows directly from the preceding proposition. For the
second formula one also applies (\ref{dd}).
$\qed$

Finally we can give the formula for the chart change on
$P^{r}_{n}X$.

\begin{prop}
In the conditions of the preceding proposition, let
$(\rho^{r}_{n}(V^{r}),(x^{i},y^{\sigma}))$
and respectively 
$(\rho^{r}_{n}(\bar{V^{r}}),(\bar{x}^{i},\bar{y}^{\sigma}))$
be the two (overlapping charts); then the change of charts on
$\rho^{r}_{n}(V^{r}) \cap \rho(\bar{V^{r}})$
is given by:

\be
\bar{y}^{\sigma}_{iI} = P^{j}_{i} D_{j} \bar{y}^{\sigma}_{I}, 
\quad |I| \leq r-1 
\label{y-y-bar}
\ee
where
$P$
is the inverse of the matrix $Q$:

\be
Q^{l}_{p} \equiv D_{p}\bar{x}^{l}, \quad P^{j}_{i} Q^{l}_{j} = \delta^{l}_{i}.
\ee
\label{change}
\end{prop}

{\bf Proof:}

We have from (\ref{invariants}) 
$$
\bar{x}^{i}_{j} \bar{y}^{\sigma}_{iI} = \bar{d}_{j} \bar{y}^{\sigma}_{I}, 
\quad |I| \leq r-1
$$
with 
$\bar{y}^{\sigma}_{I}$
functions of 
$\bar{x}^{A}_{J}.$

We will consider this relation on the overlap
$V^{r} \cap \bar{V^{r}}$
such that
$\bar{y}^{\sigma}_{I}$
can be considered as functions of
$x^{A}_{J}$
through the chart transformation formul\ae\/. Using also (\ref{dd}) one gets:
$$
\bar{x}^{i}_{j} \bar{y}^{\sigma}_{iI} = d_{j} \bar{y}^{\sigma}_{I}.
\quad |I| \leq r-1
$$

We rewrite this relation in the new coordinates
$(x^{i},y^{\sigma}_{I},x^{i}_{I})$
(see corollary \ref{new}) and also use (\ref{tilde-d}); as a result one finds
out: 

$$
z^{j}_{p} \bar{x}^{i}_{j} \bar{y}^{\sigma}_{iI} = 
\tilde{d}_{p} \bar{y}^{\sigma}_{I}, \quad |I| \leq r-1.
$$

It remains  to prove using also (\ref{tilde-d}) that

\be
Q^{j}_{p} = z^{j}_{p} \bar{x}^{i}_{j}
\label{Q}
\ee

and the change transformation formula from the statement follows.
$\qed$

We now note two other properties of the total differential operators 
$D_{i}$.
The first one follows immediately from (\ref{Dd}) and (\ref{Q}):

\be
Q^{j}_{i} \bar{D}_{j} = D_{i}.
\label{DD}
\ee

The second one is the analogue of (\ref{commutator}):

\be
\left[ D_{i}, D_{j} \right] = 0.
\ee

So the expression

$$
D_{I} \equiv \prod_{i \in I} D_{i}
$$
makes sense for every multi-index $I$.

We close this subsection with a result which will be useful later.

\begin{prop}

The following formula is true on the overlap of two charts:

\be
\Delta^{i_{1},\dots,i_{k}}_{\sigma} \bar{y}^{\nu}_{j_{1},\dots,j_{k}}
= {\cal S}^{+}_{j_{1},\dots,j_{k}} P^{i_{1}}_{j_{1}} \dots P^{i_{k}}_{j_{k}}
P^{\nu}_{\sigma}, \quad k = 1,...,r.
\label{Dy}
\ee
\end{prop}

{\bf Proof:}
It is done by recurrence. First one proves directly from the definitions that:

\be
\Delta^{i_{1},\dots,i_{k}}_{\sigma} \bar{y}^{\nu}_{j_{1},\dots,j_{k}} =
{\cal S}^{+}_{j_{1},\dots,j_{k}} P^{i_{1}}_{j_{1}} 
\Delta^{i_{2},\dots,i_{k}}_{\sigma} \bar{y}^{\nu}_{j_{2},\dots,j_{k}},
\quad k  = 2,\dots,r
\ee

and then we obtain by recurrence:

\be
\Delta^{i_{1},\dots,i_{k}}_{\sigma} \bar{y}^{\nu}_{j_{1},\dots,j_{k}} =
{\cal S}^{+}_{j_{1},\dots,j_{k}} 
P^{i_{1}}_{j_{1}}  \dots P^{i_{k-1}}_{j_{k-1}} 
\Delta^{i_{k}}_{\sigma} \bar{y}^{\nu}_{j_{k}},
\quad k = 1,\dots,r.
\ee

Finally one establishes by direct computation that

\be
\Delta^{i}_{\sigma} \bar{y}^{\nu}_{j} = P^{i}_{j} P^{\nu}_{\sigma}
\ee

and the formula from the statement follows.
$\qed$

As a corollary we have the following fact:

\begin{cor}
Let us denote by
$\Omega^{r}_{q}(PX), \quad q \geq 0$
the modulus of differential forms of order $q$ on
$P^{r}_{n}.$ 
Then the subspace 

$$
\Omega^{r}_{q,hor}(PX) \equiv \{\alpha \in \Omega^{r}_{q}(PX) |
i_{\Delta^{j_{1},\dots,j_{r}}_{\sigma}} \alpha = 0 \}
$$
is globally well defined.
\label{inner}
\end{cor}

{\bf Proof:}
One has, according to the chain rule on the overlap of two charts:

$$
\Delta^{i_{1},\dots,i_{r}}_{\sigma} =
\left( \Delta^{i_{1},\dots,i_{r}}_{\sigma} \bar{y}^{\nu}_{j_{1},\dots,j_{r}}
\right) \bar{\Delta}^{i_{1},\dots,i_{r}}_{\nu} =
P^{i_{1}}_{j_{1}} \dots P^{i_{r}}_{j_{r}} P^{\nu}_{\sigma}
\bar{\Delta}^{i_{1},\dots,i_{r}}_{\nu} 
$$

and a similar formula for the corresponding inner contractions. It follows that
the relation

$$
i_{\Delta^{j_{1},\dots,j_{r}}_{\sigma}} \alpha = 0
$$

is chart independent.
$\qed$

\subsection{Contact Forms on Grassmann Manifolds}

In this subsection we give some new material about the possibility of defining
the contact forms on the factor manifold
$P^{r}_{n}X$.
Fortunately, most of the definitions and properties from \cite{Kr4}-\cite{Kr6}
and \cite{Gr3} can be adapted to this more general situation.

By a {\it contact form} on
$P^{r}_{n}X$
we mean any form
$
\rho \in \Omega^{r}_{q}(PX)
$
verifying

\be
\left[j^{r}\gamma\right]^{*} \rho = 0
\ee
for any immersion
$\gamma: \R^{n} \rightarrow X$. 
We denote by
$\Omega^{r}_{q(c)}(PX)$
the set of contact forms of degree $q \leq n$. Here
$\left[ j^{r}\gamma\right] :\R^{n} \rightarrow P^{r}_{n}$
is given by (see def. (\ref{extension}))
$\left[ j^{r}\gamma\right] (t) \equiv \left[ j^{r}_{t}\gamma\right].$
Now, many results from \cite{Gr3} are practically unchanged. We mention some of
them.

If one considers only the 
contact forms on an open set
$\rho^{r}_{n}(V^{r}) \subset P^{r}_{n}X$
then we emphasize this by writing
$\Omega^{r}_{q(c)}(V)$.
One immediately notes that
$\Omega^{r}_{0(c)} = 0$
and that for  $q > n$ any $q$-form is contact.
It is also elementary to see that the set of all contact forms is an ideal,
denoted by
${\cal C}(\Omega^{r})$,
with respect to the operation $\wedge$. Because the operations of pull-back and
of differentiation are commuting this ideal is left invariant by exterior 
differentiation:

\be
d{\cal C}(\Omega^{r}) \subset {\cal C}(\Omega^{r}).
\label{dif-inv}
\ee

By elementary computations one finds out that, as in the case of a fibre
bundle, for any chart 
$(V,\psi)$
on $X$, every element of the set
$\Omega^{r}_{1(c)}(V)$
is a linear combination of the following expressions:

\be
\omega^{\sigma}_{j_{1},...,j_{k}} \equiv d y^{\sigma}_{j_{1},...,j_{k}} -
y^{\sigma}_{i,j_{1},...,j_{k}} d x^{i}, \quad k = 0,...,r-1
\label{o}
\ee
or, in multi-index notations

\be
\omega^{\sigma}_{J} \equiv d y^{\sigma}_{J} - y^{\sigma}_{iJ} d x^{i}, \quad
|J| \leq r-1.
\label{o'}
\ee

From the definition above it is clear that the linear subspace of the 1-forms
on 
$P^{r}_{n}$
is generated by
$d x^{i}, \quad \omega^{\sigma}_{J}, \quad (|J| \leq r-1)$
and
$dy^{\sigma}_{I}, |I| = r$.

For any smooth function on 
$\rho(V^{r})$ we have

\be
df = (d_{i}f) dx^{i} + \sum_{|J| \leq r-1} (\partial^{J}_{\sigma}f) 
\omega^{\sigma}_{J} + \sum_{|I| = r} (\partial^{I}_{\nu}f) d y^{\nu}_{I}.
\label{df'}
\ee

We also have the formula

\be
d \omega^{\sigma}_{J} = - \omega^{\sigma}_{Ji} \wedge  d x^{i}, \quad  
|J| \leq r-2.
\ee

The structure theorem from \cite{Kr6}, \cite{Gr3} stays true, i.e. any
$\rho \in \Omega^{r}_{q}(PX), \quad q = 2,...,n$
is contact {\it iff} it has the following expression in the 
associated chart:

\be
\rho = \sum_{|J|\leq r-1} \omega^{\sigma}_{J} \wedge \Phi^{J}_{\sigma} +
\sum_{|I|= r-1} d\omega^{\sigma}_{I} \wedge \Psi^{I}_{\sigma}
\label{contact-q}
\ee
where 
$\Phi^{J}_{\sigma} \in \Omega^{r}_{q-1}$
and
$\Psi^{I}_{\sigma} \in \Omega^{r}_{q-2}$
can be arbitrary forms. (We adopt the convention that 
$\Omega^{r}_{q} \equiv 0, \forall q < 0$).

We will need in the following the transformation formula relevant for change of
charts. It is to be expected that there will be some modifications of the
corresponding formula from the fibre bundle case. Namely, we have:

\begin{prop}
Let
$(V,\psi)$
and
$(\bar{V},\bar{\psi})$
two overlapping charts on $X$ and let
$(W^{r}, \Phi^{r}), \quad \Phi^{r} = (x^{i}, y^{\sigma}_{I}, x^{i}_{I})$
and
$(\bar{W}^{r}, \bar{\Phi}^{r}), \quad 
\bar{\Phi^{r}} = (\bar{x}^{i}, \bar{y}^{\sigma}_{I}, \bar{x}^{i}_{I})$
the corresponding charts on
$T^{r}_{n}X$.
Then the following formula is true on
$\rho^{r}_{n}(W^{r} \cap \bar{W}^{r})$:

\be
\bar{\omega}^{\sigma}_{I} = \sum_{|J| \leq |I|} (\Delta^{J}_{\nu}
\bar{y}^{\sigma}_{I}) \omega^{\nu}_{J} - \bar{y}^{\sigma}_{jI} (\Delta_{\nu} 
\bar{x}^{j}) \omega^{\nu}, \quad |I| \leq r-1.
\label{transf}
\ee

In particular,

\be
\bar{\omega}^{\sigma} = P^{\sigma}_{\nu} \omega^{\nu}
\label{transf0}
\ee
where we have defined:

\be
P^{\sigma}_{\nu} \equiv \Delta_{\nu} \bar{y}^{\sigma} - \bar{y}^{\sigma}_{i}
(\Delta_{\nu} \bar{x}^{i}).
\label{P}
\ee
\end{prop}

The proof goes by elementary computations. As a consequence we have:

\begin{cor}

If a $q$-form has the expression

\be
\rho  = \sum_{p+s=q-n+1} \sum_{|J_{1}|,...,|J_{p}| \leq r-1}
\sum_{|I_{1}|,...,|I_{s}|=r-1} 
\omega^{\sigma_{1}}_{J_{1}} \cdots \wedge \omega^{\sigma_{p}}_{J_{p}} \wedge
d\omega^{\nu_{1}}_{I_{1}} \cdots \wedge d\omega^{\nu_{s}}_{I_{s}} 
\wedge \Phi^{J_{1},...,J_{p},I_{1},...,I_{s}}_{\sigma_{1},...,\sigma_{p},
\nu_{1},...,\nu_{s}}
\label{str-contact-q}
\ee
is valid in one chart, then it is valid in any other chart.
\end{cor}

This corollary allows us to define for any
$q = n+1,...,N \equiv dim(J^{r}Y) = m{n+r \choose n}$  
a {\it strongly contact form} to be any
$\rho \in \Omega^{r}_{q}$ 
such that it has in one chart (thereafter in any other chart) the expression
above. For a certain uniformity of notations, we denote these forms by
$\Omega^{r}_{q(c)}$.

Now it follows that one can define the variational sequence and prove its
exactness as in the fibre bundle case.

We also mention the fact that one can define a global operator on the linear
space 
$\Omega^{r}_{q,hor}X$ 
defined at the end of the preceding subsection. In fact we have

\begin{prop}
Let 
$r \geq 1$.
Then, the operator locally defined on any differential form by:

\be
K \alpha \equiv i_{D_{j_{1}}} i_{\Delta^{j_{1},\dots,j_{r-1}}_{\sigma}} 
\left(\omega^{\sigma}_{j_{2},\dots,j_{r-1}} \wedge \alpha \right)
\ee
is globally defined on the subspace
$\Omega^{r}_{q,hor}X$.
\label{K}
\end{prop}

{\bf Proof:}
One works on the overlap of two charts and starts from the definition above
trying to transform everything into the other set of coordinates. It is quite
elementary to use corollary \ref{inner} to find

$$
K \alpha \equiv i_{\bar{D}_{j_{1}}} 
i_{\bar{\Delta}^{j_{1},\dots,j_{r-1}}_{\nu}} 
\left( P^{j_{2}}_{k_{2}} \dots P^{j_{r-1}}_{k_{r-1}} P^{\nu}_{\sigma}
\omega^{\sigma}_{j_{2},\dots,j_{r-1}} \wedge \alpha \right)
$$

Now one uses (\ref{Dy}) and the transformation formula (\ref{transf}) for the
contact forms to obtain

$$
K \alpha = \bar{K} \bar{\alpha}.
$$  
that it, $K$ is well defined globally.
$\qed$

\subsection{Morphisms of Grassmannian Manifolds}

Let 
$X_{i}, \quad i = 1,2$
be two differential manifolds and
$\phi: X_{1} \rightarrow X_{2}$
a smooth map. We define the new map
$j^{r}\phi: T^{r}_{n}X_{1} \rightarrow T^{r}_{n}X_{2}$
according to

\be
j^{r}\phi(j^{r}_{0}\gamma) \equiv j^{r}_{0} \phi \circ \gamma
\ee
for any immersion
$\gamma$. 
If 
$\gamma$
is a regular immersion, then one can see that the map
$j^{r}\phi$
maps
${\rm Imm}T^{r}_{n}(X_{1})$
into
${\rm Imm}T^{r}_{n}(X_{2})$
and so, it factorizes to a map
$J^{r}\phi: P^{r}_{n}X_{1} \rightarrow P^{r}_{n}X_{2}$
given by

\be
J^{r}\phi([j^{r}_{0}\gamma]) \equiv [j^{r}_{0} \phi \circ \gamma].
\ee

The map
$J^{r}\phi$
is called the {\it extension of order r of the map}
$\phi$.

One can show that the contact ideal behaves naturally with respect to
prolongations i.e.

\be
(J^{r}\phi)^{*}{\cal C}(\Omega^{r}(PX_{1})) \subset {\cal C}(\Omega^{r}(PX_{2})).
\label{pr-cont}
\ee

The proof follows directly from the definition of a contact form. 

If 
$\xi$
is a vector field on $X$ we define its {\it extension of order r} on
$T^{r}_{n}X$
and on
$P^{r}_{n}X$
the vector fields
$j^{r}\xi$
and 
$J^{r}\xi$
respectively given by the following formul\ae:

\be
j^{r}\xi_{j^{r}_{0}\gamma} f \equiv \left.{d \over dt} 
f \circ j^{r} e^{t\xi} (j^{r}_{0}\gamma)\right|_{t=0}
\ee
(for any smooth real function $f$ on
$T^{r}_{n}X$)
and

\be
J^{r}\xi \equiv \left.{d \over dt} J^{r} e^{t\xi}
\right|_{t=0};
\ee
here
$e^{t\xi}$
is, as usual, the flow associated to
$\xi$.

One will need the explicit formula of
$j^{r}\xi$.
If in the chart 
$(V,\psi)$
we have
\be
\xi = \xi^{A} \partial_{A} 
\label{xi}
\ee
with 
$\xi^{A}$ 
smooth function, then 
$j^{r}\xi$ 
has the following expression in the associated chart 
$(V^{r},\psi^{r})$: 

\be
j^{r}\xi = \sum_{|I| \leq r} (d_{J}\xi^{A}) \partial_{A}^{J}.
\label{pr-ev}
\ee

The proof of this fact follows by direct computation from the definition
above. We call {\it evolutions} these type of vector fields on
$T^{r}_{n}X$
and denote the set of evolutions by
${\cal E}(T^{r}_{n}X)$.

As a consequence of (\ref{pr-cont}), if 
$\xi$ 
is a vector field on 
$X$, 
then 

\be
L_{J^{r}\xi} {\cal C}(\Omega^{r}(X)) \subset {\cal C}(\Omega^{r}(X)).
\label{Lie}
\ee

Now, as in \cite{Gr3} we have the following results. Suppose that in local
coordinates we have

\be
\xi = a^{i}(x,y) {\partial \over \partial x^{i}} + b^{\sigma}(x,y)
\Delta_{\sigma} 
\ee
with $a^{i}$ and $b^{\sigma}$ smooth function; then $J^{r}\xi$ must have the
following expression in the associated chart:

\be
J^{r}\xi = a^{i}(x) {\partial \over \partial x^{i}} + \sum_{|J| \leq r} 
b^{\sigma}_{J} \Delta_{\sigma}^{J}. 
\ee
where

\be
b^{\sigma}_{J} = D_{I}(b^{\sigma} - y^{\sigma}_{j} a^{j}) + y^{\sigma}_{jI}
a^{j}, \quad |I| \leq r-1, \qquad
b^{\sigma}_{I} = D_{I}(b^{\sigma} - y^{\sigma}_{j} a^{j}), \quad |I| = r.
\ee

Finally we give the expression of the prolongation $J^{r}\phi$ where 
$\phi$ 
is a bundle morphism of the $X$. If 
$\phi$ 
has the local expression 

\be
\phi(x^{i},y^{\sigma}) = (f^{i},F^{\sigma})
\ee
then we must have in the associated chart: 

\be
J^{r}\phi(x^{i},y^{\sigma},y^{\sigma}_{j},...,y^{\sigma}_{j_{1},...,j_{r}}) =
(f^{i},F^{\sigma},F^{\sigma}_{j},...,F^{\sigma}_{j_{1},...,j_{r}})
\label{pr-transf}
\ee
where
$
F^{\sigma}_{j_{1},...,j_{k}}, \quad j_{1} \leq j_{2} \leq \cdots \leq j_{k},
\quad k = 1,...,r
$
are smooth local functions given recurringly by:

\be
F^{\sigma}_{Ji} = P^{l}_{i} D_{l}F^{\sigma}_{J}\quad |J| \leq r-1;
\ee
we also have
\be
\Delta^{I}_{\nu} F^{\sigma}_{J} = 0 \quad |I|=r.
\ee
\newpage
\section{The Lagrangian Formalism on a Grassmann Manifold}

\subsection{Euler-Lagrange Forms}

We outline a construction from \cite{An} which is the main combinatorial trick
in the study of globalisation of the Lagrangian formalism. We call any map
$P: {\cal E}(T^{r}_{n}X) \rightarrow \Omega^{s}(T^{s}_{n}X), \quad s \geq r$
covering the identity map:
$id:T^{r}_{n}X \rightarrow T^{r}_{n}X$
a {\it total differential operator}. In local coordinates such an operator has
the following expression: if 
$\xi$ 
has the local expression (\ref{xi}), then:

\be
P(\xi) = \sum_{|I| \leq r} (d_{I} \xi^{A}) P^{I}_{A} = 
\sum_{k=0}^{r} \left( d_{j_{1}} \dots d_{j_{k}} \xi^{A}\right) 
P^{j_{1},\dots,j_{k}}_{A}
\ee
where
$P^{j_{1},\dots,j_{k}}_{A}$
are differential forms on
$T^{s}_{n}X$
and
$d_{j} = d_{j}^{s}.$

Then, as in \cite{An} and \cite{Gr3} one has the following combinatorial lemma:

\begin{lemma}
In the conditions above, the following formula is true:

\be
P(\xi) = \sum_{|I| \leq r} d_{I} (\xi^{A} Q^{I}_{A})
\label{PQ}
\ee
where 

\be
Q^{I}_{A} \equiv \sum_{|J| \leq r-|I|} (-1)^{|J|} 
{|I|+|J| \choose |I|} d_{J} P^{IJ}_{A}
\label{QIA}
\ee
and one assumes that the action of a formal derivative $d_{j}$ on a form is
realized by its action on the function coefficients. Moreover, the relation
(\ref{PQ}) {\bf uniquely} determines the forms 
$Q^{I}_{\sigma}$.
\label{trick}
\end{lemma}

The proof is identical with the one presented in \cite{Gr3}. We also have

\begin{prop}
In the conditions above one has on the overlap
$V^{s} \cap \bar{V}^{s}$
the following formula:

\be
Q_{A} = (\partial_{A} \bar{x}^{B}) \bar{Q}_{B}.
\label{Q-Q-bar}
\ee

In particular, there exists a globally defined form, denoted by
$E(P)(\xi)$
with the local expression

\be
E(P)(\xi) = Q_{A} \xi^{A}.
\label{E(P)}
\ee
\end{prop}

{\bf Proof:}
From the formula (\ref{PQ}) we have

$$
P(\xi) = \xi^{A} Q_{A} +  \sum_{k=1}^{r} d_{j_{1}} \dots d_{j_{k}} 
\left( \xi^{A} Q^{j_{1},\dots,j_{k}}_{A}\right) 
$$

So, in the overlap 
$V^{s} \cap \bar{V}^{s}$
we have

$$
\xi^{A} Q_{A} - \bar{\xi}^{A} \bar{Q}_{A} =
\sum_{k=1}^{r} \left[ \bar{d}_{j_{1}} \dots \bar{d}_{j_{k}} 
\left( \bar{\xi}^{A} \bar{Q}^{j_{1},\dots,j_{k}}_{A}\right) -
d_{j_{1}} \dots d_{j_{k}} 
\left( \xi^{A} Q^{j_{1},\dots,j_{k}}_{A}\right) \right].
$$

But because of the relation (\ref{dd}) we can simplify considerably this
formula, namely we get:

$$
\xi^{A} Q_{A} - \bar{\xi}^{A} \bar{Q}_{A} =
\sum_{k=1}^{r} d_{j_{1}} \dots d_{j_{k}} 
\left( \bar{\xi}^{A} \bar{Q}^{j_{1},\dots,j_{k}}_{A} -
\xi^{A} Q^{j_{1},\dots,j_{k}}_{A}\right)
$$

Now one proves that both sides are zero as in \cite{An}, \cite{Gr3} making
use of Stokes theorem.
$\qed$

The operator $E(P)$ defined by (\ref{E(P)}) is called the {\it Euler operator}
associated to the total differential operator $P$; it has the local expression:

\be
E(P)(\xi) = \xi^{A} E_{A}(P)
\ee
where

\be
E_{A}(P) = \sum_{|I|=0}^{r} (-1)^{|I|} d_{I} P^{I}_{A}.
\ee

Now one takes 
${\cal L} \in {\cal F}(T^{r}_{n})$
and constructs the total differential operator
$P_{\cal L}$
according to:

\be
P_{\cal L}(\xi) \equiv L_{pr(\xi)} {\cal L}.
\ee

Lemma \ref{trick} can be applied and immediately gives the following local
formula:

\be
P_{\cal L}(\xi) = \sum_{|I|=0}^{r} d_{I} \left( \xi^{A}
{\cal E}^{I}_{A}({\cal L})\right)
\ee
where

\be
{\cal E}^{I}_{A}(L) \equiv \sum_{|J|\leq r-|I|} (-1)^{|J|} 
{|I|+|J| \choose |I|} d_{J} \partial^{IJ}_{A} {\cal L}
\label{Lie-Euler}
\ee
are the so-called {\it Lie-Euler operators}; the Euler operator associated
to 
$P_{\cal L}$ 
has the following expression:

\be
E(P_{\cal L}) = \xi^{A} {\cal E}_{A}(L)
\ee
where

\be
{\cal E}_{A}(L) \equiv \sum_{|J|\leq r} (-1)^{|J|} 
d_{J} \partial^{J}_{A} {\cal L}
\label{EL1}
\ee
are the {\it Euler-Lagrange expressions} associated to 
${\cal L}$.

The proposition above leads to

\begin{prop}
If
${\cal L} \in {\cal F}(T^{r}_{n})$,
then there exists a globally defined $1$-form, denoted by
${\cal E}({\cal L})$
such that we have in the chart 
$V^{s}, \quad s \geq 2r$:

\be
{\cal E}({\cal L}) = {\cal E}_{A}({\cal L}) dx^{A}.
\label{EL2}
\ee
\end{prop}

{\bf Proof:}
By construction
$$
Q_{A} = {\cal E}_{A}({\cal L})
$$

Now, one has from (\ref{Q-Q-bar}) :

\be
{\cal E}_{A}({\cal L}) = (\partial_{A} \bar{x}^{B}) 
\bar{\cal E}_{B}(\bar{\cal L}).
\label{EE}
\ee

Combining with the transformation property

$$
d\bar{x}^{A} = (\partial_{B} \bar{x}^{A}) dx^{B}
$$
and obtains that the formula (\ref{EL2}) has a global meaning.
$\qed$

One calls this form the {\it Euler-Lagrange form} associated to 
${\cal L}$.

All the properties of this form listed in \cite{Gr3} are true in this case
also. We insist only on the so-called {\it product rule} for the Euler-Lagrange
expressions \cite{Al} which will be repeatedly used in the following.

\begin{prop}
If $f$ and $g$ are smooth functions on
$V^{r}$
then one has in
$V^{s}, \quad s \geq 2r$
the following formula:

\be
{\cal E}^{I}_{A}(fg) = \sum_{|J| \leq r-|I|} {|I| + |J| \choose |J|}
\left[ (d_{J} f) {\cal E}^{IJ}_{A}(g) + (d_{J} g) {\cal E}^{IJ}_{A}(f)\right],
\quad |I| \leq r.
\ee 
\label{product}
\end{prop}

The proof goes by direct computation, directly from the definition of the
Lie-Euler operators combined with Leibnitz rule of differentiation of a
product. 

We now come to the main definition. A smooth real function 
${\cal L}$ 
on
${\rm Imm}T^{r}_{n}$
is called a {\it homogeneous Lagrangian} if it verifies the relation:

\be
{\cal L}(x\cdot a) = det(a) {\cal L}(x), \quad \forall a \in L^{r}_{n};
\label{homogeneous}
\ee
here by
${\rm det}(a)$
we mean
${\rm det}(a^{i}_{j})$.

Such an object induces on the factor manifold
$P^{r}_{n}X$
an non-homogeneous object.

\begin{prop}

Let
${\cal L}$
be a homogeneous Lagrangian. Then for every chart
$(W^{r},\Phi^{r})$
on
${\rm Imm}T^{r}_{n}X$
there exists a smooth real function on 
$\rho^{r}_{n}(W^{r})$
such that:

\be
{\cal L} = {\rm det}({\bf x}) L \circ \rho^{r}_{n}
\label{LL}
\ee
\end{prop}
and conversely, if 
${\cal L}$
is locally defined by this relation, then it verifies (\ref{homogeneous}).

{\bf Proof:}
One chooses in (\ref{homogeneous})
$a = {\bf x}$
(see (\ref{x})).
$\qed$

The function $L$ is called the {\it non-homogeneous (local) Lagrangian}
associated to 
${\cal L}$.

As a consequence of the connection (\ref{LL}) we have

\begin{prop}

Let
${\cal L}$
a homogeneous Lagrangian and
$(V,\psi), \quad (\bar{V},\bar{\psi})$
two overlapping charts on $X$. We consider on the associated charts
$(W^{r},\Phi^{r})$
and
$(\bar{W}^{r},\bar{\Phi}^{r})$
the corresponding non-homogeneous Lagrangians
$L$ and respectively
$\bar L$.
Then we have on the overlap
$\rho^{r}_{n}(W^{r} \cap \bar{W}^{r})$
the following formula:

\be
L = {\cal J} \bar L
\label{ll}
\ee
where

\be
{\cal J} \equiv {\rm det}(Q) = {\rm det}(D_{i}\bar{x}^{j}).
\ee
\label{l1}
\end{prop}

{\bf Proof:}
One writes (\ref{LL}) for both charts and gets

$$
{\cal L} = {\rm det}({\bf x}) L \circ \rho^{r}_{n} = 
{\rm det}(\bar{\bf x}) \bar{L} \circ \rho^{r}_{n};
$$
as a consequence

$$
L \circ \rho^{r}_{n} = {\rm det}(\bar{\bf x} {\bf z}) \bar{L} \circ
\rho^{r}_{n}.  $$

One now uses the relation (\ref{Q}) and obtains the relation from the
statement.  
$\qed$

As consequence we have

\begin{thm}

One can globally define the equivalence class

$$
[\lambda] \in \Omega^{r}_{n}(PX) /\Omega^{r}_{n(c)}(PX) 
$$
such that the local expression of $\lambda$ is

\be
\lambda = L \theta_{0};
\ee
here, as usual

\be
\theta_{0} \equiv dx^{1} \wedge \dots \wedge dx^{n}.
\ee
\label{l2}
\end{thm}

{\bf Proof:}
One proves immediately that on the overlap of the associated charts from the
proposition above one has:

\be
\bar{\theta}_{0} = {\cal J} \theta_{0} + {\rm contact \quad terms}.
\ee

This result must be combined with (\ref{ll}) to obtain on
$\rho^{r}_{n}(W^{r} \cap \bar{W}^{r})$:

$$
\lambda - \bar{\lambda} \in \Omega^{r}_{n(c)}(PX);
$$
this proves the theorem.
$\qed$

It is natural to ask what is the connection between the Euler-Lagrange
expression of the homogeneous and the corresponding non-homogeneous Lagrangian.
The answer is contained in:

\begin{thm}

Suppose 
${\cal L}$
is a homogeneous Lagrangian defined on
${\rm Imm}T^{r}_{n}X$
and
$L$
is the associated non-homogeneous Lagrangian. Then the following relations are 
valid on the chart 
$W^{s}, \quad s \geq 2r$:

\be
{\cal E}_{\sigma}({\cal L}) = 
{\rm det}({\bf x}) E_{\sigma}(L) \circ \rho^{s}_{n},
\label{e1}
\ee

\be
{\cal E}^{j_{1},\dots,j_{k}}_{\sigma}({\cal L}) = (-1)^{k} {\rm det}({\bf x}) 
\sum_{|I| \geq k} (-1)^{|I|} \sum\limits_{(I_{1},\dots,I_{k})} 
{\cal S}^{+}_{j_{1},\dots,j_{k}} z^{j_{1}}_{I_{1}} \dots z^{j_{k}}_{I_{k}}
E^{I}_{\sigma}(T) \circ \rho^{s}_{n},
\label{e2}
\ee

and

\be
{\cal E}_{p}({\cal L}) = 
- {\rm det}({\bf x}) y^{\sigma}_{p} E_{\sigma}(T) \circ \rho^{s}_{n},
\label{e3}
\ee

\be
{\cal E}^{j}_{p}({\cal T}) = {\rm det}({\bf x}) 
\left[ z^{j}_{p} L + y^{\sigma}_{p} \sum_{I \leq r} z^{j}_{I}
E^{I}_{\sigma}(L) \right] \circ \rho^{s}_{n},
\label{e4}
\ee

\be
{\cal E}^{j_{1},\dots,j_{k}}_{p}({\cal L}) = (-1)^{k+1} {\rm det}({\bf x}) 
y^{\sigma}_{p} \sum_{|I|\geq k} (-1)^{|I|} 
\sum\limits_{(I_{1},\dots,I_{k})} 
{\cal S}^{+}_{j_{1},\dots,j_{k}} z^{j_{1}}_{I_{1}} \dots z^{j_{k}}_{I_{k}}
E^{I}_{\sigma}(L) \circ \rho^{s}_{n}.
\label{e5}
\ee
\label{e}
\end{thm}

{\bf Proof:}

(i) As a general strategy of the proof, we will try to transform the expression
of the total differential operator 
$P_{\cal L}(\xi)$
in terms of $L$; we have by definition

$$
P_{\cal L}(\xi) = \sum_{|I| \leq r} (d_{I} \xi^{A}) \partial^{I}_{A}  
\left[ {\rm det}({\bf x}) L \circ \rho\right].
$$

Using the chain rule one obtains rather easily from here:

\begin{eqnarray}
P_{\cal L}(\xi) = {\rm det}({\bf x}) \sum_{k=0}^{r} 
(d_{j_{1}} \dots d_{j_{k}} \xi^{\sigma}) \sum_{l=k}^{r}
(\Delta^{i_{1},\dots,i_{l}}_{\nu} L) \circ \rho 
\quad (\partial^{j_{1},\dots,j_{k}}_{\sigma} y^{\nu}_{i_{1},\dots,i_{l}}) +
\nonumber \\
{\rm det}({\bf x}) \sum_{k=0}^{r} 
(d_{j_{1}} \dots d_{j_{k}} \xi^{p}) \sum_{l=k}^{r}
(\Delta^{i_{1},\dots,i_{l}}_{\nu} L) \circ \rho 
\quad (\partial^{j_{1},\dots,j_{k}}_{p} y^{\nu}_{i_{1},\dots,i_{l}}) +
\nonumber \\
{\rm det}({\bf x}) \xi^{p} (\partial_{p} L) \circ \rho +
(d_{j}\xi^{p}) [\partial^{j}_{p} {\rm det}({\bf x})] L \circ \rho 
\label{PL}
\end{eqnarray}

We will now consider that the functions
$\xi^{A}$
depend only of the variables
$(x^{i},y^{\sigma})$
i.e. there exist the smooth real functions
$\Xi^{A}$
on
$W^{r}$
such that

\be
\xi^{A} = \Xi^{A} \circ \rho.
\ee

(ii) To compute further the expression above one starts from (\ref{xy}) and
firstly proves directly:

\be
\partial^{j_{1},\dots,j_{k}}_{\sigma} y^{\nu}_{I} = 
\delta^{\nu}_{\sigma} {\cal S}^{+}_{j_{1},\dots,j_{k}}
\sum\limits_{(I_{1},\dots,I_{k})} z^{j_{1}}_{I_{1}} \cdots z^{j_{k}}_{I_{k}},
\quad 1 \leq k \leq |I| \leq r.
\label{dy1}
\ee

Next one proves

\begin{lemma}

If we have 
$1 \leq k \leq |I| \leq r$
then one has:

\be
\partial^{j_{1},\dots,j_{k}}_{p} y^{\sigma}_{I} = 
- {\cal S}^{+}_{j_{1},\dots,j_{k}}
\sum\limits_{(I_{0},\dots,I_{k})} z^{j_{1}}_{I_{1}} \cdots z^{j_{k}}_{I_{k}}
y^{\sigma}_{pI_{0}}.
\label{dy2}
\ee

Here we understand that the subsets
$I_{1},\dots,I_{k}$
cannot be the emptyset; on the contrary, it is allowed to have
$I_{0} = \emptyset.$
\end{lemma}

{\bf Proof:}
By induction on 
$|I|$
starting from
$|I| = k$.
For this smallest possible value one uses (\ref{xy}). Then one supposes that 
the formula from the statement is valid for
$k \leq |I| < r$,
uses the defining recurrence relation (\ref{invariants}) for the invariants and
establishes the relation for
$iI$. The cases
$k = 1$
and $k > 1$
must be treated separately.
$\nabla$

Another auxiliary result is contained in the well known result:

\begin{lemma}

The following formul\ae\/ are true:

\be
\partial^{i}_{j} {\rm det}({\bf x}) = z^{i}_{j} {\rm det}({\bf x})
\ee

and

\be
d_{i} [z^{i}_{j} {\rm det}({\bf x})] = 0.
\ee
\end{lemma}

{\bf Proof:}
The first result is quite general i.e. valid for any invertible matrix and it
is proved directly from the definition of the determinant. The second formula
is a corollary of the first.
$\nabla$

One must use these results together with (\ref{dD}); we have from (\ref{PL}) 
after permuting the two summation signs: 

\begin{eqnarray}
P_{\cal L}(\xi) = {\rm det}({\bf x}) 
[\xi^{\sigma} (\Delta_{\sigma} L) \circ \rho +
\sum_{|I|=1}^{r}  (\Delta^{I}_{\sigma} L) \circ \rho \quad \sum_{k=1}^{I}
\sum\limits_{(I_{1},\dots,I_{k})} z^{j_{1}}_{I_{1}} \cdots z^{j_{k}}_{I_{k}}
( d_{j_{1}} \dots d_{j_{k}} \xi^{\sigma}) -
\nonumber \\
\sum_{|I_{0}| \leq r} \sum_{|I|=1}^{r-|I_{0}|} {|I|+|I_{0}| \choose |I|}
y^{\sigma}_{pI_{0}} (\Delta^{II_{0}}_{\sigma} L) \circ \rho \quad 
\sum_{k=1}^{I} \sum\limits_{(I_{1},\dots,I_{k})} 
z^{j_{1}}_{I_{1}} \cdots z^{j_{k}}_{I_{k}} 
(d_{j_{1}} \dots d_{j_{k}} \xi^{p}) +
\nonumber \\
\xi^{p} (\partial_{p} L) \circ \rho - \xi^{p} (D_{j} L) \circ \rho ] +
d_{j}[ \xi^{p} z^{j}_{p} {\rm det}({\bf x}) L \circ \rho]. 
\label{PL1}
\end{eqnarray}

(iii) One can proceed further one must generalise the formula (\ref{dD}).
We have the following formula

\be
d_{I}(f \circ \rho) = \sum_{k=1}^{I} \sum\limits_{(I_{1},\dots,I_{k})}
x^{j_{1}}_{I_{1}} \cdots x^{j_{k}}_{I_{k}}
(D_{j_{1}} \dots D_{j_{k}} f) \circ \rho
\ee
which can be proved by induction on
$|I|$.

This formula can be ``inverted" rather easily and we get:

\be
(D_{I}f) \circ \rho = \sum_{k=1}^{I} \sum\limits_{(I_{1},\dots,I_{k})}
z^{j_{1}}_{I_{1}} \cdots z^{j_{k}}_{I_{k}}
d_{j_{1}} \dots d_{j_{k}} (f \circ \rho).
\ee

The expression (\ref{PL1}) can be considerably simplified to

\be
P_{\cal L}(\xi) = {\rm det}({\bf x}) [\sum_{|I|\leq 1}
(D_{I} \Xi^{\sigma}) P^{I}_{\sigma} - (D_{I} \Xi^{p}) P^{I}_{p}] \circ \rho +
d_{j}[ \xi^{p} z^{j}_{p} {\rm det}({\bf x}) L \circ \rho] 
\label{PL2}
\ee
where we have defined

$$
P^{I}_{\sigma} \equiv \Delta^{I}_{\sigma} L; \quad
P^{I}_{p} \equiv \sum_{|J| \leq r-|I|} {|I|+|J| \choose |J|} y^{\sigma}_{pJ}
\Delta^{IJ}_{\sigma}.
$$

(iv) Now it this the time to apply lemma \ref{trick} and to obtain in this way:

\be
P_{\cal L}(\xi) = {\rm det}({\bf x}) \sum_{|I|\leq r}
[D_{I} (\Xi^{\sigma} Q^{I}_{\sigma} - \Xi^{p} Q^{I}_{p})] \circ \rho +
d_{j}[ \xi^{p} z^{j}_{p} {\rm det}({\bf x}) L \circ \rho] 
\label{PL3}
\ee
where one gets by some computations the following formul\ae\/ for the
expressions 
$Q^{I}_{A}$:

\be
Q^{I}_{\sigma} = \sum_{|J| \leq r-|I|} (-1)^{|J|} {|I|+|J| \choose |J|} 
D_{J} P^{IJ}_{\sigma} = E^{I}_{\sigma}(L)
\ee

and

\be
Q^{I}_{p} = \sum_{|J| \leq r-|I|} (-1)^{|J|} {|I|+|J| \choose |J|} 
D_{J} P^{IJ}_{p} = y^{\sigma}_{p} E^{I}_{\sigma}(L).
\ee

(v) We want to compare the expression (\ref{PL3}) with the right hand side of
the formula from lemma \ref{trick}. To do this we need one more combinatorial
result valid for any smooth real function on the chart
$W^{r}$:

\be
{\rm det}({\bf x}) (D_{I} f) \circ \rho = \sum_{k=1}^{I}
(-1)^{r-|I|} d_{j_{1}} \cdots d_{j_{k}} \left[ {\rm det}({\bf x}) 
\sum\limits_{(I_{1},\dots,I_{k})} z^{j_{1}}_{I_{1}} \cdots z^{j_{k}}_{I_{k}}
(f \circ \rho) \right], \quad |I| \leq r.
\ee

One proves this formula by induction on 
$|I|$
and so the final expression for the total differential operator is

\begin{eqnarray}
P_{\cal L}(\xi) = 
{\rm det}({\bf x}) (\Xi^{\sigma} Q_{\sigma} - \Xi^{p} Q^{I}_{p}) \circ \rho +
d_{j}[ \xi^{p} z^{j}_{p} {\rm det}({\bf x}) L \circ \rho] + \nonumber \\
\sum_{k=1}^{r} (-1)^{k} d_{j_{1}} \cdots d_{j_{k}} \left[ {\rm det}({\bf x}) 
\sum_{|I|\geq k} (-1)^{|I|} 
\sum\limits_{(I_{1},\dots,I_{k})} z^{j_{1}}_{I_{1}} \cdots z^{j_{k}}_{I_{k}}
(\Xi^{\sigma} Q^{I}_{\sigma} - \Xi^{p} Q^{I}_{p}) \circ \rho \right]
\label{PL4}
\end{eqnarray}

If one uses the unicity statement from lemma \ref{trick} one obtains the
desired formul\ae.
$\qed$

Immediate consequences of the preceding theorem are

\begin{cor}

If 
${\cal L}$
is a homogeneous Lagrangian, then the following relations are true on the 
manifold
$T^{s}_{n}X, \quad s \geq 2r:$

\be
{\cal E}_{A}({\cal L}) (x\cdot a) = {\rm det}(a) {\cal E}_{A}({\cal L}) (x),
\quad \forall a \in L^{r}_{n}.
\ee

In the conditions of the above theorem we have:

\be
{\cal E}_{A}({\cal L}) \equiv 0 \Longleftrightarrow E_{\sigma}(L) \equiv 0.
\ee
\label{hom-e}
\end{cor}

We also have the analogue of proposition \ref{l1}:

\begin{prop}

In the condition of the preceding theorem let us consider two overlapping
charts 
$(V,\psi), \quad (\bar{V},\bar{\psi})$.  
Then one has on the overlap of the associated charts: 
$\rho^{s}_{n}(W^{s} \cap \bar{W}^{s})$ 
the following relation:

\be
E_{\sigma}(L) = {\cal J} P^{\nu}_{\sigma} \bar{E}_{\nu}(\bar{L})
\ee
where the matrix
$P^{\nu}_{\sigma}$
has been defined by the formula (\ref{P}).
\label{el1}
\end{prop}

{\bf Proof:}

We start from the transformation formula (\ref{EE}) for
$A \rightarrow \sigma$:

$$
{\cal E}_{\sigma}({\cal L}) = 
(\Delta_{\sigma} \bar{x}^{i}) \bar{\cal E}_{i}(\bar{\cal L}) + 
(\Delta_{\sigma} \bar{y}^{\nu}) \bar{\cal E}_{\nu}(\bar{\cal L})
$$
and substitute (\ref{e1}) and (\ref{e3}). Using (\ref{Q}) we obtain by 
elementary computations the relation from the statement.
$\qed$

\begin{rem}
For a different proof of this result see \cite{Ol}.
\end{rem}

Now we have the analogue of theorem \ref{l2}:

\begin{thm}

If ${\cal L}$ is a homogeneous Lagrangian on
$T^{r}_{n}$,
then one can globally define the equivalence class

$$
[E(L)] \in \Omega^{s}_{n+1}(PX) /\Omega^{s}_{n+1(c)}(PX)
$$
on
$P^{s}_{n}, \quad s\geq 2r$
such that the local expression of 
$E(L)$ is

\be
E(L) = E_{\sigma}(L) \omega^{\sigma} \wedge \theta_{0}.
\ee
\label{el2}
\end{thm}

{\bf Proof:} Follows the lines of theorem \ref{l2} and it is elementary.
$\qed$

We also note the following property:

\begin{prop}

If ${\cal L}$ is a homogeneous Lagrangian on
$T^{r}_{n}X$,
then the corresponding Euler-Lagrange form verifies on
$T^{s}_{n}X, \quad s \geq 2r$
the following identity:

\be
(j^{s}_{0}\gamma)^{*} {\cal E}({\cal L}) = 0, \quad 
\forall \gamma \in {\rm Imm}T^{s}_{n}.
\ee
\label{cond}
\end{prop}

{\bf Proof:}
By direct computation we get

$$
(j^{s}_{0}\gamma)^{*} {\cal E}({\cal L}) = 
\left( {\cal E}_{A}({\cal L}) x^{A}_{j} \right) \circ j^{s}_{0}\gamma 
\quad dt^{j}.
$$

Now one uses (\ref{e1}) and (\ref{e3}) to prove that the expression in the
bracket is identically zero.
$\qed$

We close this subsection with some remarks.

\begin{rem}
If ${\cal L}$ is a homogeneous Lagrangian, one can expect some homogeneity
property for the total differential operator associated to it. Indeed, one has
for an arbitrary Lagrangian 

\be
P_{\cal L} \circ \phi_{a}(\xi) = P_{\cal L}(\xi) \circ \phi_{a}
\ee
where
$\phi_{a}$
denotes the right action of the differential group 
$L^{r}_{n}$.

As a consequence, one has for a homogeneous Lagrangian

\be
P_{\cal L}(\xi) \circ \phi_{a} = {\rm det}(a) P_{\cal L}(\xi), \quad
\forall a \in L^{r}_{n}.
\ee
\end{rem}

\subsection{Differential Equations on Grassmann Manifolds}

An element 
${\cal T} \in \Omega^{s}_{n+1}(X)$
is called a {\it differential equation} on
$T^{s}_{n}X$. 
In the chart
$(V^{s},\psi^{s})$
the differential equation ${\cal T}$ has the following local expression:

\be
{\cal T} = {\cal T}_{A} dx^{A}.
\label{T}
\ee

It is clear that the Euler-Lagrange form defined by (\ref{EL2}) is a
differential equation. A differential equation ${\cal T}$ is called {\it
variational} if there exists a Lagrangian ${\cal L}$ on
$T^{r}_{n}X, \quad s \geq 2r$
such that we have
${\cal T} = {\cal E}({\cal L})$. 
If the function ${\cal L}$ is only locally defined, then such a differential
equation is called {\it locally variational} (or, of the {\it Euler-Lagrange
type}). 

If 
$\gamma:\R^{n} \rightarrow X$ 
is a immersion, then on says that the differential equation ${\cal T}$
{\it verifies the differential equation} {\it iff} we have

\be
(j^{s}_{0}\gamma)^{*} i_{Z} {\cal T} = 0
\ee
for any vector field $Z$ on $J^{s}_{n}X$. In local coordinates we have on 
$V^{s}$:

\be
{\cal T}_{A} \circ j^{s}_{0}\gamma = 0 \quad (A = 1,...,N).
\label{eq-diff}
\ee

Guided by corollary \ref{hom-e} and prop \ref{cond} we also introduce the
following definition.  We say that ${\cal T}$  is a {\it homogeneous 
differential equation} on
${\rm Imm}T^{s}_{n}X$ 
if it verifies the following conditions:

\be
(\phi_{a})^{*} {\cal T} = {\rm det}(a) {\cal T}, \quad \forall a \in L^{s}_{n}
\ee

and

\be
(j^{s}_{0}\gamma)^{*} {\cal T} = 0, \quad 
\forall \gamma \in {\rm Imm}T^{s}_{n}.
\ee

Then we have the following result which can be proved by elementary
computations suggested by the similar results obtained for a differential
equation of the Euler-Lagrange type.

\begin{thm}

Let ${\cal T}$ be a homogeneous differential equation on
${\rm Imm}T^{s}_{n}X$.
Then there exist some local smooth real functions
$T_{\sigma}$
in every chart
$\rho^{s}_{n}(W^{s})$
such that one has:

\be
{\cal T}_{\sigma} = {\rm det}({\bf x}) T_{\sigma} \circ \rho, \qquad
{\cal T}_{i} = - {\rm det}({\bf x}) y^{\sigma}_{i} T_{\sigma} \circ \rho.
\ee

As a consequence, if 
$(V,\psi), \quad (\bar{V},\bar{\psi})$
are two overlapping charts on $X$, then one has on the intersection
$\rho^{s}_{n}(W^{s} \cap \bar{W}^{s})$
the following transformation formula:

\be
T_{\sigma} = {\cal J} P^{\nu}_{\sigma} \bar{T}_{\nu}
\ee

and the class

$$
[T] \in \Omega^{s}_{n+1}(PX) /\Omega^{s}_{n+1(c)}(PX)
$$
can be properly defined such that we have locally

\be
T = T_{\sigma} \omega^{\sigma} \wedge \theta_{0}.
\ee
\label{homogeneous-eq}
\end{thm}

One calls $T$ the {\it associated (local) non-homogenous differential 
equation}.

We now prove the existence of the (globally) defined Helmholtz-Sonin form
associated to a differential equation. By analogy with \cite{Gr3} we have the 
following result

\begin{thm}
Let ${\cal T}$ be a differential equation on
$T^{s}_{n}X$
with the local form
given by (\ref{T}). We define the following expressions in any chart 
$V^{t}, \quad t > 2s$:

\be
{\cal H}^{J}_{AB} \equiv \partial^{J}_{B} {\cal T}_{A} - (-1)^{|J|}
{\cal E}^{J}_{A}({\cal T}_{B}), \quad |J| \leq s.
\label{HS-local}
\ee

Then there exists a globally defined $2$-form, denoted by 
${\cal H}({\cal T})$ 
such that in any chart $V^{t}$ we have:

\be
{\cal H}({\cal T}) = \sum_{|J| \leq s} {\cal H}^{J}_{AB} dx^{B}_{J} \wedge
dx^{A}.
\label{HS}
\ee
\end{thm}

{\bf Proof:}
We sketch briefly the argument from \cite{Gr3}. Let $\xi$ be a vector field on 
$X$; we define a (global) $1$-form 
${\cal H}_{\xi}({\cal T})$ 
according to:

\be
{\cal H}_{\xi}({\cal T}) \equiv L_{pr(\xi)} {\cal T} - 
{\cal E}\left( i_{pr(\xi)} {\cal T}\right)
\label{H-xi}
\ee
and the following local expression is obtained:

\be
{\cal H}_{\xi}({\cal T}) = \sum_{|I| \leq s} (d_{I} \xi^{B}) 
{\cal H}^{I}_{AB} dx^{A}.
\ee

The transformation formula at a change of charts for the expressions 
$d_{I} \xi^{B}$ 
is:

\be
\bar{d}_{I} \bar{\xi}^{A} = \sum_{|J| \leq |I|} \left(\partial^{J}_{B} 
\bar{x}^{A}_{I} \right) (d_{J} \xi^{B}), \quad |I| = 0,...,s.
\label{d-xi}
\ee
 
Using the transformation formula (\ref{d-xi}) one can obtain the
transformation formula for the expressions 
${\cal H}^{I}_{AB}$: 
one has in
the overlap 
$V^{t} \cap \bar{V}^{t}, \quad t \geq 2s$:

\be
{\cal H}^{J}_{DB} = \sum_{|I| \geq |J|} 
\left( \partial^{J}_{B} \bar{x}^{A}_{I} \right) 
\left( \partial_{D} \bar{x}^{C} \right) \bar{\cal H}^{I}_{CA}.
\label{HH}
\ee

This transformation formula leads now to the fact that
${\cal H}({\cal T})$ 
has an invariant meaning.
$\qed$
 
${\cal H}({\cal T})$ 
is called the {\it Helmholtz-Sonin form} associated to ${\cal T}$ and 
${\cal H}^{I}_{AB}$ 
are the {\it Helmholtz-Sonin expressions} associated to ${\cal T}$. 

A well-known corollary of the theorem above is:

\begin{cor}
The differential equation ${\cal T}$ is locally variational {\it iff} 
${\cal H}({\cal T}) = 0$
iff

\be
{\cal H}^{I}_{AB} = 0, \quad \forall A,B = 1,\dots N, \quad \forall |I| \leq r.
\label{HSeq}
\ee
\end{cor}

The proof is identical with the one presented in \cite{Gr3}. The preceding
equations are called the {\it Helmholtz-Sonin equations}.

As in the preceding subsection, if ${\cal T}$ is a homogeneous differential
equation, we have a very precise connection between the Helmholtz-Sonin
expressions of ${\cal T}$ and of $T$ from theorem \ref{homogeneous-eq}.

\begin{thm}

Suppose 
${\cal T}$
is a homogeneous equation defined on
${\rm Imm}T^{s}_{n}X$
and
$T_{\sigma}$
the components of the associated non-homogeneous equation. Then the following
relations are valid on the chart 
$W^{t}, \quad t \geq 2s$:

\be
{\cal H}_{\sigma\nu}({\cal T}) = 
{\rm det}({\bf x}) H_{\sigma\nu}(T) \circ \rho^{t}_{n},
\label{hs1}
\ee

\be
{\cal H}^{j_{1},\dots,j_{k}}_{\sigma\nu}({\cal T}) = {\rm det}({\bf x}) 
\sum_{|I| \geq k} \sum\limits_{(I_{1},\dots,I_{k})} 
{\cal S}^{+}_{j_{1},\dots,j_{k}} z^{j_{1}}_{I_{1}} \dots z^{j_{k}}_{I_{k}}
H^{I}_{\sigma\nu}(T) \circ \rho^{t}_{n}, \quad k = 1,\dots,s.
\label{hs2}
\ee

\be
{\cal H}_{\sigma p}({\cal T}) = 
- {\rm det}({\bf x}) \sum_{I_{0} \leq s} y^{\nu}_{pI_{0}} 
H^{I_{0}}_{\sigma\nu}(T) \circ \rho^{t}_{n},
\label{hs3}
\ee

\be
{\cal H}^{j_{1},\dots,j_{k}}_{\sigma p}({\cal T}) = - {\rm det}({\bf x}) 
\sum_{I_{0} \leq s} \sum_{|I|=k}^{s-|I_{0}|} {|I|+|I_{0}| \choose |I|} 
\sum\limits_{(I_{1},\dots,I_{k})} 
{\cal S}^{+}_{j_{1},\dots,j_{k}} z^{j_{1}}_{I_{1}} \dots z^{j_{k}}_{I_{k}}
y^{\nu}_{pI_{0}} H^{II_{0}}_{\sigma\nu}(T) \circ \rho^{t}_{n},
\quad k = 1,\dots,s,
\label{hs4}
\ee

\be
{\cal H}_{p\nu}({\cal T}) = 
- {\rm det}({\bf x}) y^{\sigma}_{p} H_{\sigma\nu}(T) \circ \rho^{t}_{n},
\label{hs5}
\ee

\be
{\cal H}^{j_{1},\dots,j_{k}}_{p\nu}({\cal T}) = - {\rm det}({\bf x}) 
y^{\sigma}_{p} \sum_{I_{0} \leq s} \sum_{|I| \geq k} 
\sum\limits_{(I_{1},\dots,I_{k})} 
{\cal S}^{+}_{j_{1},\dots,j_{k}} z^{j_{1}}_{I_{1}} \dots z^{j_{k}}_{I_{k}}
H^{I}_{\sigma\nu}(T) \circ \rho^{t}_{n}, \quad k = 1,\dots,s
\label{hs6}
\ee

and

\be
{\cal H}_{pq}({\cal T}) = 
{\rm det}({\bf x}) y^{\sigma}_{p} \sum_{I_{0} \leq s} y^{\nu}_{qI_{0}} 
H^{I_{0}}_{\sigma\nu}(T) \circ \rho^{t}_{n},
\label{hs7}
\ee

\be
{\cal H}^{j_{1},\dots,j_{k}}_{pq}({\cal T}) = {\rm det}({\bf x}) 
y^{\sigma}_{p} \sum_{I_{0} \leq s} \sum_{|I|=k}^{s-|I_{0}|}
{|I|+|I_{0}| \choose |I|} \sum\limits_{(I_{1},\dots,I_{k})} 
{\cal S}^{+}_{j_{1},\dots,j_{k}} z^{j_{1}}_{I_{1}} \dots z^{j_{k}}_{I_{k}}
y^{\nu}_{qI_{0}} H^{II_{0}}_{\sigma\nu}(T) \circ \rho^{t}_{n}, 
\quad k = 1,\dots,s.
\label{hs8}
\ee
\label{hs}
\end{thm}

The proof is tedious but elementary. One must use the formul\ae\/ derived in
theorem \ref{e} combined with the derivation property from proposition 
\ref{product} to prove case by case the formul\ae\/ from the statement.
Occasionally, one must study separately the cases
$k = 0, 1$
and
$k > 1$.

In the conditions of the above theorem we have:

\begin{cor}

\be
{\cal H}^{I}_{AB}({\cal T}) \equiv 0 \Longleftrightarrow 
H_{\sigma\nu}(T) \equiv 0.
\label{HSHS}
\ee
\end{cor}

As it can be expected we have the analogues of propositions \ref{l1} and
\ref{el1}:

\begin{prop}

In the condition of the preceding theorem let us consider two overlapping 
charts
$(V,\psi), \quad (\bar{V},\bar{\psi})$.
Then one has on the overlap of the associated charts:
$\rho^{t}_{n}(W^{t} \cap \bar{W}^{t})$
the following relations:

\be
H_{\sigma\nu}(T) = {\cal J} P^{\alpha}_{\sigma} \sum_{|J| \leq s}
\left[ \left(\Delta_{\nu} 
\bar{y}^{\beta}_{J}\right) - P^{i}_{j}\left( D_{i} \bar{y}^{\beta}_{J}\right)
\left( \Delta_{\nu} \bar{x}^{j}\right) \right] 
\bar{H}^{J}_{\alpha\beta}(\bar{T})
\label{tr1}
\ee

\be
H^{I}_{\sigma\nu}(T) = {\cal J} P^{\alpha}_{\sigma} \sum_{|J| \geq |I|}
\left(\Delta_{\nu} \bar{y}^{\beta}_{J}\right) 
\bar{H}^{J}_{\alpha\beta}(\bar{T}), \quad \forall I \not= \emptyset.
\label{tr2}
\ee
\label{hhss1}
\end{prop}

{\bf Proof:}
(i) It is convenient to introduce the expression
$g_{j_{1},\dots,j_{k}}, \quad k = 0,\dots,r$
completely symmetric in all indices (with the convention
$g_{\emptyset} = g$)
and to use the (\ref{hs1}) and (\ref{hs2}) to obtain:

\begin{eqnarray}
{\rm det}({\bf x}) \sum_{|I| \leq s} (g\cdot {\bf z})_{I} 
H^{I}_{\sigma\nu}(T) \rho = \nonumber \\
\sum_{k=0}^{s} g_{j_{1},\dots,j_{k}} 
{\cal H}^{j_{1},\dots,j_{k}}_{\sigma\nu}({\cal T}) = \nonumber \\
\sum_{|I| \leq s} g_{I} \sum_{|J| \geq |I|} 
\left( \partial^{J}_{\sigma} \bar{x}^{A}_{I} \right) 
\left( \partial_{\nu} \bar{x}^{C} \right) \bar{\cal H}^{I}_{CA}(\bar{\cal T})
\nonumber
\end{eqnarray}
where use have been made of the transformation formula (\ref{HH}) for
$B \rightarrow \sigma, \quad D \rightarrow \nu$.

If we make here
$g \rightarrow g\cdot {\bf x}$
we obtain equivalently after elementary prelucrations of the formula above:

\begin{eqnarray}
{\rm det}({\bf x}) H_{\sigma\nu}(T) \circ \rho  = 
\sum_{|I| \leq s}
\sum_{|J| \geq |I|} [(\partial_{\nu} \bar{x}^{\beta}_{I}) 
(\partial_{\sigma} \bar{x}^{\alpha}) 
\bar{\cal H}^{I}_{\alpha\beta}(\bar{\cal T}) +
(\partial_{\nu} \bar{x}^{\beta}_{I}) 
(\partial_{\sigma} \bar{x}^{p}) \bar{\cal H}^{I}_{p\beta}(\bar{\cal T}) + 
\nonumber \\
(\partial_{\nu} \bar{x}^{q}_{I}) 
(\partial_{\sigma} \bar{x}^{\alpha}) 
\bar{\cal H}^{I}_{\alpha q}(\bar{\cal T})  + 
(\partial_{\nu} \bar{x}^{q}_{I}) 
(\partial_{\sigma} \bar{x}^{p}) \bar{\cal H}^{I}_{pq}(\bar{\cal T})] 
\label{hhh1}
\end{eqnarray}

and for 
$k \geq 1$:

\begin{eqnarray}
{\rm det}({\bf x}) H^{j_{1},\dots,j_{k}}_{\sigma\nu}(T) \circ \rho  =
\sum_{|J| \geq k} \sum\limits_{(J_{1},\dots,J_{k})} 
{\cal S}^{+}_{j_{1},\dots,j_{k}} x^{j_{1}}_{J_{1}} \dots x^{j_{k}}_{J_{k}}
\sum_{|J| \geq |I|} [(\partial^{J}_{\nu} \bar{x}^{\beta}_{I}) 
(\partial_{\sigma} \bar{x}^{\alpha}) 
\bar{\cal H}^{I}_{\alpha\beta}(\bar{\cal T}) +
\nonumber \\
(\partial^{J}_{\nu} \bar{x}^{\beta}_{I}) 
(\partial_{\sigma} \bar{x}^{p}) \bar{\cal H}^{I}_{p\beta}(\bar{\cal T}) + 
(\partial^{J}_{\nu} \bar{x}^{q}_{I}) 
(\partial_{\sigma} \bar{x}^{\alpha}) 
\bar{\cal H}^{I}_{\alpha q}(\bar{\cal T}) + 
(\partial^{J}_{\nu} \bar{x}^{q}_{I}) 
(\partial_{\sigma} \bar{x}^{p}) \bar{\cal H}^{I}_{pq}(\bar{\cal T})]
\label{hhh2}
\end{eqnarray}

The second relation can be considerably simplified if one uses (\ref{yx}), more
precisely the consequence

\be
\Delta^{j_{1},\dots,j_{k}}_{\nu} x_{I}^{\mu} = \delta^{\mu}_{\nu}
\sum\limits_{(I_{1},\dots,I_{k})}
{\cal S}^{+}_{j_{1},\dots,j_{k}}
x_{I_{1}}^{j_{1}} \dots x_{I_{k}}^{j_{k}};
\ee
the chain rule gives from here

\be
\Delta^{j_{1},\dots,j_{k}}_{\nu} \bar{x}_{I}^{A} = \sum_{|J|=k}^{|I|}
(\partial_{\nu}^{J} \bar{x}^{A}_{I})
\sum\limits_{(J_{1},\dots,J_{k})}
{\cal S}^{+}_{j_{1},\dots,j_{k}} x_{J_{1}}^{j_{1}} \dots x_{J_{k}}^{j_{k}}.
\ee

So, the relation (\ref{hhh2}) becomes:

\begin{eqnarray}
{\rm det}({\bf x}) H^{j_{1},\dots,j_{k}}_{\sigma\nu}(T) \circ \rho  = 
\sum_{|I| \geq k} 
\{(\Delta^{j_{1},\dots,j_{k}}_{\nu} \bar{x}^{\beta}_{I}) 
[(\partial_{\sigma} \bar{x}^{\alpha}) 
\bar{\cal H}^{I}_{\alpha\beta}(\bar{\cal T}) + \nonumber \\
(\partial_{\sigma} \bar{x}^{p}_{I}) \bar{\cal H}^{I}_{p\beta}(\bar{\cal T})] +
(\Delta^{j_{1},\dots,j_{k}}_{\nu} \bar{x}^{q}_{I}) 
[(\partial_{\sigma} \bar{x}^{\alpha}) \bar{\cal H}^{I}_{\alpha q}(\bar{\cal T})
+ (\partial_{\sigma} \bar{x}^{p}) \bar{\cal H}^{I}_{pq}(\bar{\cal T}) ]\}
\label{hhh}
\end{eqnarray}

If we compare with (\ref{hhh1}) we see that the preceding relation stays true
for 
$ k = 0$
also.

Now we use again the theorem above in the right hand side of the relation
just derived and obtains after some computations (using the relations
(\ref{dy1}) and (\ref{dy2}) and the chain rule) the relations from the
statement of the theorem.
$\qed$

Now we have the analogue of theorems \ref{l2} and \ref{el2}

\begin{thm}

If ${\cal T}$ is a homogeneous differential equation on
${\rm Imm}T^{s}_{n}$,
then one can globally define the equivalence class

$$
[H(T)] \in \Omega^{t}_{n+2}(PX) /\Omega^{t}_{n+2(c)}(PX) 
$$
on
$P^{t}_{n}, \quad t\geq 2s$
such that the local expression of 
$H(T)$ is

\be
H(T) = \sum_{|I| \leq s} H^{I}_{\sigma\nu}(L) \omega^{\nu}_{I} \wedge 
\omega^{\sigma} \wedge \theta_{0}.
\ee
\label{hhss2}
\end{thm}

{\bf Proof:} Follows the lines of theorem \ref{el2} and it is elementary.
$\qed$

As a consequence of the theorems \ref{l2}, \ref{el2} and \ref{hhss2} we can
apply the exactness of the variational sequence and obtain that the expressions
of an arbitrary variationally trivial Lagrangian and of a locally variational
differential equation are coinciding with the expressions derived in
\cite{Gr4}. 
\newpage
\section{Lagrangian Formalism on Second Order Grassmann Bundles}

\subsection{The Second Order Grassmann Bundle}

Here we particularize the results obtained in the preceding sections for the 
case
$r = 2$.

The coordinates on
$T^{2}_{n}X$
are
$(x^{A},x^{A}_{j},x^{A}_{ij})$
and with the help of the derivative operators (see (\ref{partial}))

\be
\partial_{A} \equiv {\partial \over \partial x^{A}}, \quad
\partial_{A}^{j} \equiv {\partial \over \partial x^{A}_{j}}, \quad
\partial_{A}^{ij} \equiv 
\cases{{\partial \over \partial x^{A}_{ij}}, & {\rm for} $i = j$ \cr
{1 \over 2} {\partial \over \partial x^{A}_{ij}}, & {\rm for} $i \not= j$ \cr}
\label{partial2}
\ee
we have for any smooth function $f$ (see (\ref{df})):

\be
df = (\partial^{A} f) dx^{A} + 
(\partial^{A}_{i} f) dx^{A}_{i} +
(\partial^{A}_{ij} f) dx^{A}_{ij}
\label{df2}
\ee

We have the following formul\ae\/ (see (\ref{derivation})):

\be
\partial_{A} x^{B} = \delta^{B}_{A}, \quad
\partial^{i}_{A} x^{B}_{j} = \delta^{B}_{A} \delta^{i}_{j}, \quad
\partial^{ij}_{A} x^{B}_{lm} = {1 \over 2} \delta^{B}_{A}
(\delta^{i}_{l} \delta^{j}_{m} + \delta^{i}_{m} \delta^{j}_{l})
\label{derivation2}
\ee
and the other derivatives are zero.
The formal derivatives (see (\ref{formal})) are in this case:

\be
d^{r}_{i} \equiv x^{A} \partial_{A} + x^{A}_{ij} \partial^{j}_{A}
\label{formal2}
\ee
and from here we immediately have (see (\ref{d-x})):

\be
d_{i} x^{A} = x^{A}_{i}, \quad 
d_{i} x^{A}_{j} = x^{A}_{ij}.  
\label{d-x2}
\ee

The formul\ae\/ for the induces change of charts (see (\ref{change-chart})) are
in this case:

\be
F_{i}^{A} = x_{i}^{B} \partial_{B} F^{A}, \quad 
F_{i_{1},i_{2}}^{A} = x_{i_{1},i_{2}}^{B} \partial_{B} F^{A} +
x_{i_{1}}^{B_{1}} x_{i_{2}}^{B_{2}} \partial_{B_{1}} \partial_{B_{2}} F^{A}
\label{change-chart2}
\ee

The elements of the differential group are of the form

\be
a = (a^{j}_{i},a^{j}_{i_{1},i_{2}}), \quad {\rm det}(a^{j}_{i}) \not= 0
\ee
with the composition law (see (\ref{composition})):

\be
(a \cdot b)^{k}_{i} = b^{k}_{j} a^{j}_{i}, \quad
(a \cdot b)^{k}_{i_{1},i_{2}} = b^{j}_{i_{1},i_{2}} a^{k}_{j} + 
b^{j_{1}}_{i_{1}} b^{j_{2}}_{i_{2}} a^{k}_{j_{1},j_{2}}
\label{composition2}
\ee
and the inverse element given by

\be
a^{-1} = ((a^{-1})^{j}_{i},- (a^{-1})^{j}_{k} (a^{-1})^{j_{1}}_{i_{1}}
(a^{-1})^{j_{2}}_{i_{2}} a^{k}_{j_{1},j_{2}}).
\ee

The action of this group on 
$T^{2}_{n}X$
is (see (\ref{action})):

\be
(a \cdot x)^{A} = x^{A}, \quad
(a \cdot x)^{A}_{i} = a^{j}_{i} x^{A}_{j}, \quad
(a \cdot x)^{A}_{i_{1},i_{2}} = a^{j}_{i_{1},i_{2}} x^{A}_{j} +
a^{j_{1}}_{i_{1}} a^{j_{2}}_{i_{2}} x^{A}_{j_{1},j_{2}}.
\ee

The expressions for the invariants of this action are (see (\ref{xy})):

\be
y^{\sigma} = x^{\sigma}, \quad
y_{i}^{\sigma} = z_{i}^{j} x_{j}^{\sigma}, \quad
y_{i_{1},i_{2}}^{\sigma} = z_{i_{1}}^{j_{1}} z_{i_{2}}^{j_{2}} 
(x_{j_{1},j_{2}}^{\sigma} - z^{k}_{p} y^{\sigma}_{k} x_{j_{1},j_{2}}^{p})
\label{xy2}
\ee
and they are, together with
$x^{i}$,
local coordinates on
$P^{2}_{n}X.$

The inverse of these formul\ae\/ are (see (\ref{yx})):

\be
x^{\sigma} = y^{\sigma}, \quad
x_{i}^{\sigma} = x_{i}^{j} y_{j}^{\sigma}, \quad
x_{i_{1},i_{2}}^{\sigma} = x_{i_{1}}^{j_{1}} x_{i_{2}}^{j_{2}} 
y_{j_{1},j_{2}}^{\sigma} + y^{\sigma}_{k} x_{i_{1},i_{2}}^{k}
\label{yx2}
\ee

On the factor manifold
$P^{2}_{n}X$
one introduces the derivatives operators (see (\ref{delta}) and (\ref{D})):

\be
\Delta_{\sigma} \equiv {\partial \over \partial y^{\sigma}}, \quad
\Delta_{\sigma}^{j} \equiv {\partial \over \partial y_{j}^{\sigma}}, \quad
\Delta_{\sigma}^{ij} \equiv 
\cases{{\partial \over \partial y^{\sigma}_{ij}}, & {\rm for} $i = j$ \cr
{1 \over 2} {\partial \over \partial y^{\sigma}_{ij}}, & {\rm for} $i \not= j$ 
\cr}
\label{delta2}
\ee

and 

\be
D_{i} \equiv {\partial \over \partial x^{i}} + 
y^{\sigma}_{i} \Delta_{\sigma} + y^{\sigma}_{ij} \Delta^{j}_{\sigma}.
\label{D2}
\ee

The formula for the change of charts on
$P^{2}_{n}X$
is (see (\ref{y-y-bar})):

\be
\bar{y}^{\sigma}_{i} = P^{j}_{i} D_{j} \bar{y}^{\sigma}, \quad
\bar{y}^{\sigma}_{i_{1},i_{2}} = P^{j_{1}}_{i_{1}} P^{j_{2}}_{i_{2}} 
\left[ - P^{m}_{l} (D_{j_{1}} D_{j_{2}} \bar{x}^{l}) (D_{m} \bar{y}^{\sigma}) +
D_{j_{1}} D_{j_{2}} \bar{y}^{\sigma}\right].
\label{y-y-bar2}
\ee

Finally, the expressions for the contact forms are (see (\ref{o})): 

\be
\omega^{\sigma} \equiv d y^{\sigma} - y^{\sigma}_{i} d x^{i}, \quad 
\omega^{\sigma}_{j} \equiv d y^{\sigma}_{j} - y^{\sigma}_{ij} d x^{i}
\label{o2}
\ee

and their transformation for a change of charts is (see (\ref{transf}) and
(\ref{transf0})):

\be
\bar{\omega}^{\sigma} = P^{\sigma}_{\nu} \omega^{\nu}, \quad
\bar{\omega}^{\sigma}_{i} = P^{\sigma}_{\nu} P^{j}_{i} 
\bar{\omega}^{\nu}_{j} +  
\left[\Delta_{\nu} \bar{y}^{\sigma}_{i} - P^{l}_{k}(\Delta_{\nu} \bar{x}^{k}) 
(D_{l} \bar{y}^{\sigma}_{i})\right] \omega^{\nu}.
\label{transf2}
\ee

\subsection{Lagrangian Formalism}

Here we give a different approach to the Lagrangian formalism based on a
certain 
$(n+1)$-form
defined on the Grassmann manifold
$P^{2}_{n}$
\cite{GP}. The description of the formalism will be slightly different a some
new material will appear. 

As in \cite{GP}, \cite{Gr2}, we base our formalism on the operator $K$ (see
prop. \ref{K}) which in our case is defined on
$\Omega^{2}_{q,hor}(PX)$
and is given by

\be
K \alpha \equiv i_{D_{j}} i_{\Delta^{j}_{\sigma}} 
\left(\omega^{\sigma} \wedge \alpha \right).
\label{K2}
\ee

We define the space of {\it Lagrange-Souriau forms} according to

\be
\Omega^{2}_{LS} \equiv \{ \alpha \in \Omega^{2}_{n+1,hor}(PX) | d\alpha = 0, 
\quad K\alpha = 0, \quad i_{V_{1}} i_{V_{2}} \alpha = 0,
\quad \forall V_{i} \in {\rm Vert}(X) \}.
\label{closed}
\ee
where
${\rm Vert}(X)$
is the space of vertical vector fields on 
$P^{2}_{n}$
with respect to the projection
$\rho^{2,1}_{n}$.

By definition, a {\it Lagrangian system} on
$P^{2}_{n}$
is a couple
$(E,\alpha)$
where $E$ is a open sub-bundle of
$P^{2}_{n}$
and
$\alpha$
is a Lagrange-Souriau form.

If 
$\gamma \in {\rm Imm}T^{2}_{n}$
we say that it {\it verifies the Euler-Lagrange equations iff} 

\be
[j^{2}_{0}\gamma]^{*} \alpha = 0.
\label{eqEL2}
\ee
 
It is easy to see that the local expression of a Lagrange-Souriau form is

\begin{eqnarray}
\alpha = \sum_{k=0}^{n} {1 \over k!}
F^{i_{0},\dots,i_{k}}_{\sigma_{0},\dots,\sigma_{k}}
\omega^{\sigma_{0}}_{i_{0}} \wedge \omega^{\sigma_{1}} \wedge \cdots 
\omega^{\sigma_{k}} \wedge \theta_{i_{1},\dots,i_{k}} + \nonumber \\
\sum_{k=0}^{n} {1 \over (k+1)!}
E^{i_{1},\dots,i_{k}}_{\sigma_{0},\dots,\sigma_{k}}
\omega^{\sigma_{0}} \wedge \omega^{\sigma_{1}} \wedge \cdots 
\omega^{\sigma_{k}} \wedge \theta_{i_{1},\dots,i_{k}}
\label{alpha}
\end{eqnarray}
where we have defined:

\be
\theta_{i_{1},\dots,i_{k}} \equiv {n \choose k} \varepsilon_{i_{1},\dots,i_{n}}
dx^{i_{k+1}} \wedge \dots \wedge dx^{i_{n}}, \quad k = 0,...,n.
\ee

We can admit, without loosing generality, some (anti)-symmetry properties.

\be
F^{i_{0},i_{P(1)},\dots,i_{P(k)}}_{\sigma_{0},\sigma_{Q(1)},\dots,
\sigma_{Q(k)}} = (-1)^{|P|+|Q|}
F^{i_{0},i_{1},\dots,i_{k}}_{\sigma_{0},\sigma_{1},\dots,\sigma_{k}},
\quad \forall P, Q \in {\cal P}_{k}
\label{anti-F}
\ee
and

\be
E^{i_{P(1)},\dots,i_{P(k)}}_{\sigma_{Q(0)},\dots,\sigma_{Q(k)}} = 
(-1)^{|P|+|Q|}
E^{i_{1},\dots,i_{k}}_{\sigma_{0},\sigma_{1},\dots,\sigma_{k}},
\quad \forall P \in {\cal P}_{k}, \forall Q \in {\cal P}_{k+1}.
\label{anti-E}
\ee

The condition 
$$
K \alpha = 0
$$
appearing in the definition (\ref{closed}) of a Lagrange-Souriau form, has the
following local form (see \cite{Gr2}):

\be
{\cal S}^{-}_{i_{0},\dots,i_{k}} {\cal S}^{-}_{\sigma_{0},\dots,\sigma_{k}}
F^{i_{0},\dots,i_{k}}_{\sigma_{0},\dots,\sigma_{k}} = 0, \quad k = 1,\dots,n
\label{LS}
\ee
where 
${\cal S}^{\pm}_{\sigma_{0},\dots,\sigma_{k}}$
are defined similarly to (\ref{sa}).

The local form of the Euler-Lagrange equation (\ref{eqEL2}) is simply:

\be
E_{\sigma} \circ [j^{2}_{0}\gamma] = 0
\label{local-EL}
\ee

\begin{rem}
For the case
$n = 1$,
the $2$-form $\alpha$ above appears in \cite{Su} and the condition (\ref{LS})
is investigated in \cite{Kl} and \cite{Ho}.
\end{rem}

The justification of the terminology for (\ref{local-EL}) is contained in the 
following result:

\begin{prop}

The expressions
$E_{\sigma}$
verify the Helmholtz equations (\ref{HSHS}).
\end{prop}

{\bf Proof:}
One writes in detail the closedness condition 
$d\alpha = 0$
and find out, in particular, the following equations:

\be
d_{j} E^{j}_{\sigma_{0},\sigma_{1}} + \partial_{\sigma_{0}} E_{\sigma_{1}} -
\partial_{\sigma_{1}} E_{\sigma_{0}} = 0,
\label{HSa}
\ee

\be
d_{l} F^{jl}_{\nu,\sigma} - \partial_{\sigma} F^{j}_{\nu} +
\partial^{j}_{\nu} E_{\sigma} + E^{j}_{\nu,\sigma} = 0,
\label{HSb}
\ee

\be
\partial^{jk}_{\nu} E_{\sigma} + {1\over 2} \left(F^{jk}_{\nu,\sigma} +
F^{kj}_{\nu,\sigma}\right) = 0,
\label{HSc}
\ee

\be
F^{jl}_{\sigma,\nu} - F^{lj}_{\nu,\sigma} = 
\partial^{l}_{\nu} F^{j}_{\sigma} - \partial^{j}_{\sigma} F^{l}_{\nu}.
\label{HSd}
\ee

From (\ref{LS}) we get, in particular:

\be
F^{j}_{\sigma} = 0
\label{LSa}
\ee
and

\be
F^{jl}_{\nu,\sigma} - F^{jl}_{\sigma,\nu} - F^{lj}_{\nu,\sigma} +
F^{lj}_{\sigma,\nu} = 0.
\label{LSb}
\ee

We use (\ref{LSa}) in (\ref{HSb}) and get:

$$
d_{l} F^{jl}_{\nu,\sigma} + \partial^{j}_{\nu} E_{\sigma} + 
E^{j}_{\nu,\sigma} = 0,
$$

If we substitute the last term of this relation into (\ref{HSa}) and use
(\ref{HSc}) we get:

\be
\partial_{\sigma_{0}} E_{\sigma_{1}} - \partial_{\sigma_{1}} E_{\sigma_{0}} = 
d_{j} \partial^{j}_{\sigma_{0}} E_{\sigma_{1}} -
d_{j} d_{l} \partial^{jl}_{\sigma_{0}} E_{\sigma_{1}}.
\label{h1}
\ee

Next, we take the symmetric part in 
$\nu, \sigma$
of (\ref{HSb}) and obtain:

$$
\partial^{j}_{\nu} E_{\sigma} + \partial^{j}_{\sigma} E_{\nu} =
- d_{l} \left( F^{jl}_{\nu,\sigma} + F^{jl}_{\sigma,\nu}\right).
$$

One uses here (\ref{LSb}) and next (\ref{HSd}) + (\ref{LSa}) to get:

\be
\partial^{j}_{\nu} E_{\sigma} + \partial^{j}_{\sigma} E_{\nu} =
2 d_{l} \partial^{jl}_{\sigma} E_{\nu}.
\label{h2}
\ee

Finally, the antisymmetric part of (\ref{HSc}) in
$\sigma, \nu$
is:

\be
\partial^{jk}_{\nu} E_{\sigma} = \partial^{jk}_{\sigma} E_{\nu}.
\label{h3}
\ee

The equations (\ref{h1}), (\ref{h2}) and (\ref{h3}) are the Helmholtz-Sonin
equations for the expressions
$E_{\sigma}.$
$\qed$

We now mention a result derived in \cite{GP}:

\begin{prop}
There exists in every chart a local $n$-form 
$\beta$
on
$P^{1}_{n}X$
having the local expression

\be
\beta = \sum_{k=0}^{n} {1 \over k!}
L^{i_{1},\dots,i_{k}}_{\sigma_{1},\dots,\sigma_{k}}
\omega^{\sigma_{1}} \wedge \cdots \omega^{\sigma_{k}} \wedge
\theta_{i_{1},\dots,i_{k}}
\label{beta}
\ee
where

\be
L^{i_{1},\dots,i_{k}}_{\sigma_{1},\dots,\sigma_{k}} = 
{\cal S}^{-}_{i_{1},\dots,i_{k}} {\cal S}^{-}_{\sigma_{1},\dots,\sigma_{k}}
\partial^{i_{1}}_{\sigma_{1}} \dots \partial^{i_{k}}_{\sigma_{k}} L, 
\quad k = 0,\dots, n
\label{beta-L}
\ee

such that:

\be
\alpha = d(\rho^{2,1}_{n})^{*}\beta.
\label{alpha-beta}
\ee
\end{prop}

\begin{rem}
The form $\beta$ is a generalization of the Poincar\'e-Cartan form \cite{Ca},
\cite{Po}. It had appeared in the literature in \cite{Kr1}, \cite{Be1}, 
\cite{Be2}, \cite{Ru2}. For other generalisations of the Poincar\'e-Cartan form
see \cite{Ru1}, \cite{Ki}, \cite{GS}, \cite{Ga}, \cite{Go}, \cite{Sa} and
\cite{Kr2} where the notion of Lepage is introduced for such generalisations.
\end{rem}

As a consequence we can express the coefficients of the form $\alpha$ given by
(\ref{alpha}) in terms of the smooth function $L$ \cite{Gr2}:

\be
F^{i_{0},\dots,i_{k}}_{\sigma_{0},\dots,\sigma_{k}} =
\partial^{i_{0}}_{\sigma_{0}} 
L^{i_{1},\dots,i_{k}}_{\sigma_{1},\dots,\sigma_{k}} -
L^{i_{0},\dots,i_{k}}_{\sigma_{0},\dots,\sigma_{k}}, \quad k = 0,\dots,n
\label{F-sigma}
\ee
and

\be
E^{i_{1},\dots,i_{k}}_{\sigma_{0},\dots,\sigma_{k}} =
{\cal S}^{-}_{i_{0},\dots,i_{k}} \partial_{\sigma_{0}}
L^{i_{1},\dots,i_{k}}_{\sigma_{1},\dots,\sigma_{k}} -
d_{i_{0}} L^{i_{0},\dots,i_{k}}_{\sigma_{0},\dots,\sigma_{k}},
\quad k = 0,\dots,n.
\label{E-sigma}
\ee

In particular, for 
$k = 0$
the preceding formula is

\be
E_{\sigma} = \partial_{\sigma} L - d_{i} \partial^{i}_{\sigma} L
\ee
i.e. the Euler-Lagrange operator. This is another justification of the
terminology for the equations (\ref{eqEL2}) (and (\ref{local-EL}).) This shows
that the expressions $E^{i_{1},\dots,i_{k}}_{\sigma_{0},\dots,\sigma_{k}},
\quad k = 0,\dots,n$ are also some generalizations of the Euler-Lagrange
expressions, however, different from the Lie-Euler expressions introduced in
\cite{Al} and given by (\ref{Lie-Euler}).

The coefficients of $\alpha$ verify some recurrence relations:

\be
E^{i_{1},\dots,i_{k}}_{\sigma_{0},\dots,\sigma_{k}} = -
{\cal S}^{-}_{i_{1},\dots,i_{k}} {\cal S}^{-}_{\sigma_{0},\dots,\sigma_{k}} 
\left( \partial^{i_{1}}_{\sigma_{0}} + d_{l} \partial^{li_{1}}_{\sigma_{0}}
\right) E^{i_{2},\dots,i_{k}}_{\sigma_{1},\dots,\sigma_{k}}, 
\quad k = 1,\dots,n
\ee
and

\be
F^{i_{0},\dots,i_{k}}_{\sigma_{0},\dots,\sigma_{k}} = {\rm const} \times
{\cal S}^{-}_{\sigma_{0},\dots,\sigma_{k}} {\cal S}^{-}_{i_{0},\dots,i_{k}}
\partial^{i_{0}}_{\sigma_{0}}
F^{i_{1},\dots,i_{k}}_{\sigma_{1},\dots,\sigma_{k}}, \quad k = 2,\dots,n.
\ee

{\bf Proof:}
By exploiting the conditions
$d\alpha = 0$
and 
$K\alpha = 0$
written in local coordinates. Alternatively, one computes the right hand sides
of these formul\ae\/ using (\ref{F-sigma}) and (\ref{E-sigma}) and obtains the
left hand sides.  
$\qed$

These recurrence relations have an important consequence:

\begin{cor}
In the conditions above we have

\be
\alpha \equiv 0 \Longleftrightarrow E_{\sigma} \equiv 0.
\ee
\end{cor}

We close this subsection giving the connection between the form $\alpha$
introduced here and the Lagrange-Souriau form $\sigma$ introduced in
\cite{Gr2}, \cite{GP}. We have

\begin{prop}
Let $\alpha$ be a Lagrange-Souriau form. Then there exists a 
$(n+1)$-form $\sigma$ on
$P^{1}_{n}$
such that

\be
\alpha = (\rho^{2,1}_{n})^{*}\sigma.
\label{alpha-sigma}
\ee
\end{prop}

{\bf Proof:}
One explicitates the condition
$$
d \alpha = 0
$$
appearing in the definition of a Lagrange-Souriau form (\ref{closed}) and 
finds out in particular:

\be
\partial^{jk}_{\nu} F^{i_{0},\dots,i_{k}}_{\sigma_{0},\dots,\sigma_{k}} = 0,
\quad k = 0,\dots,n
\ee

\be
\partial^{jk}_{\nu} E^{i_{1},\dots,i_{k}}_{\sigma_{0},\dots,\sigma_{k}} +
{1 \over 2} \left[ F^{j,k,i_{0},\dots,i_{k}}_{\nu,\sigma_{0},\dots,\sigma_{k}} 
+ (j \leftrightarrow k) \right] = 0, \quad k = 0,\dots,n.
\ee

It follows that the generic expression of the coefficients
$E^{i_{1},\dots,i_{k}}_{\sigma_{0},\dots,\sigma_{k}}$
is:

\be
E^{i_{1},\dots,i_{k}}_{\sigma_{0},\dots,\sigma_{k}} = 
G^{i_{1},\dots,i_{k}}_{\sigma_{0},\dots,\sigma_{k}} -
y^{\nu}_{jk} F^{j,k,i_{1},\dots,i_{k}}_{\nu,\sigma_{0},\dots,\sigma_{k}}
\label{EG}
\ee
where
$G^{i_{1},\dots,i_{k}}_{\sigma_{0},\dots,\sigma_{k}}$
have the antisymmetry property (\ref{anti-E}) and verify

\be
\partial^{jk}_{\nu} G{i_{0},\dots,i_{k}}_{\sigma_{0},\dots,\sigma_{k}} = 0,
\quad k = 0,\dots,n.
\ee

Substituting (\ref{EG}) into the expression (\ref{alpha}) one obtains
(\ref{alpha-sigma}) with

\begin{eqnarray}
\sigma = \sum_{k=0}^{n} {1 \over k!}
F^{i_{0},\dots,i_{k}}_{\sigma_{0},\dots,\sigma_{k}}
dy^{\sigma_{0}}_{i_{0}} \wedge \omega^{\sigma_{1}} \wedge \cdots 
\omega^{\sigma_{k}} \wedge \theta_{i_{1},\dots,i_{k}} + \nonumber \\
\sum_{k=0}^{n} {1 \over (k+1)!}
G^{i_{1},\dots,i_{k}}_{\sigma_{0},\dots,\sigma_{k}}
\omega^{\sigma_{0}} \wedge \omega^{\sigma_{1}} \wedge \cdots 
\omega^{\sigma_{k}} \wedge \theta_{i_{1},\dots,i_{k}}.
\label{sigma}
\end{eqnarray}

This finishes the proof.
$\qed$

It is plausible that the formalism presented in this section can be extended 
to Grassmann manifolds of arbitrary order $r > 2$. Some steps in this
direction are contained in \cite{Kr5}.

\vskip 1cm

{\bf Acknowledgments:} 
The conjecture that it is possible to define the Lagrange, Euler-Lagrange and
Helmholtz-Sonin classes in this more general setting of non-fibred manifolds was
suggested to the author by professor Demeter Krupka.


\newpage
\begin{thebibliography}{99}

\bibitem{Al}
S. J. Aldersley,
``{\it Higher Euler Operators and some of their Applications}",
Journ. Math. Phys. {\bf 20} (1979) 522-531
 
\bibitem{An}
I. M. Anderson,
``{\it The Variational Bicomplex}",
Utah State Univ. preprint, 1989, (Academic Press, Boston, to appear)

\bibitem{AD}
I. M. Anderson, T. Duchamp,
``{\it On the Existence of Global Variational Principles}",
American Journ. Math. {\bf 102} (1980) 781-868

\bibitem{Be1} 
D. E. Betounes, 
``{\it Extensions of the Classical Cartan Form}", 
Phys. Rev. {\bf D 29} (1984) 599-606

\bibitem{Be2} 
D. E. Betounes, 
``{\it Differential Geometric Aspects of the Cartan Form: Symmetry Theory}", 
J. Math. Phys. {\bf 28} (1987) 2347-2353

\bibitem{Ca} 
E. Cartan, 
``{\it Le\c cons sur les Invariants Integraux}", 
Hermann, 1922.

\bibitem{D}
J. Dieudonn\'e,
``{\it \'El\'ements d'Analyse 3},
Gauthier-Villars, Paris, 1970

\bibitem{Ga} 
P. L. Garcia, 
``{\it The Poincar\'e-Cartan Invariant in the Calculus of Variations}", 
Symp. Math. {\bf XIV} (1974) 219-246

\bibitem{GS} 
H.Goldschmits, S.Sternberg: 
``{\it The Hamilton-Cartan Formalism in the Calculus of Variations}", 
Ann. Inst. Fourier {\bf 23} (1973) 203-267

\bibitem{Go} 
M. J. Gotay, 
``{\it A Multisymplectic Framework for Classical Field Theory and the 
Calculus of Variations. I. Covariant Hamiltonian  Formalism}", 
in ``Mechanics, Analysis and Geometry: 200 Years after Lagrange", 
M. Francaviglia and D. D. Holms, eds., North-Holland, Amsterdam, 1990, 
pp. 203-235

\bibitem{Gr2}
D. R. Grigore,
``{\it A Generalized Lagrangian Formalism in Particle Mechanics and Classical
Field Theory}",
Fortschr. der Phys. {\bf 41} (1993) 569-617

\bibitem{Gr3}
D. R. Grigore,
``{\it Variational Sequence on Finite Jet Bundle Extensions and the Lagrangian 
Formalism}",
dg-ga/9702016, submitted for publication

\bibitem{Gr4} 
D. R. Grigore,
``{\it Variationally Trivial Lagrangians and Locally Variational Differential 
Equations of Arbitrary Order}", 
submitted for publication, 

\bibitem{GK} 
D. R. Grigore, D. Krupka,
``{\it Invariants of Velocities and Higher Order Grassmann Bundles}",
dg-ga/9708013
to appear in Journ. Geom. Phys.

\bibitem{GP} 
D. R. Grigore, O. T. Popp, 
``{\it On the Lagrange-Souriau Form in Classical Field Theory}", 
to appear in Mathematica Bohemica

\bibitem{Ho} 
P. Horv\`athy,
``{\it Variational Formalism for Spinning Particles}",
Journ. Math. Phys. {\bf 20} (1979) 49-52

\bibitem{Kl} 
J. Klein, 
``{\it Espaces Variationels et M\'ecaniqu}e", 
Ann. Inst. Fourier {\bf 12} (1962) 1-124

\bibitem{Ki} 
I. Kijowski, 
``{\it A Finite-Dimensional Canonical Formalism in Classical Field Theory}", 
Comm. Math. Phys. {\bf 30} (1973) 99-128

\bibitem{Kr1} 
D. Krupka, 
``{\it A Map Associated to the Lepagean Forms of the Calculus 
of Variations in Fibered Manifolds}", 
Czech. Math. Journ. {\bf 27} (1977) 114-118

\bibitem{Kr2}
D. Krupka,
``{\it Lepagean Forms in Higher Order Variational Theory}",
in Proceedings of the IUTAM-ISIMM Symposium on ``Modern Developments in
Analytical Mechanics", Turin, 1982, Atti della Academia delle Scienze di
Torino, Suppl. al Vol. {\bf 117} (1983) 198-238

\bibitem{Kr4}
D. Krupka,
``{\it Variational Sequence on Finite Order Jet Spaces}",
in Proceedings of the Conference ``Differential Geometry and its Applications",
August, 1989, World Scientific, Singapore, 1990, pp. 236-254

\bibitem{Kr5}
D. Krupka,
``{\it Topics in the Calculus of Variation: Finite Order Variational
Sequences}", 
in ``Differential Geometry and its Applications", proceedings Conf. Opava,
1992, Open Univ. Press, pp. 437-495

\bibitem{Kr6}
D. Krupka,
``{\it The Contact Ideal}",
Diff. Geometry and its Applications {\bf 5} (1995) 257-276

\bibitem{Ol}
P. J. Olver,
``{\it Applications of Lie Groups to Differential Equations}",
Springer, 1986

\bibitem{Po} 
H. Poincar\'e, 
``{\it Le\c cons sur les M\'ethodes Nouvelles de la  M\'ecanique C\'eleste}", 
Gauthier-Villars, Paris, 1892.

\bibitem{Ru1} 
H. Rund, 
``{\it A Cartan Form for the Field Theory of Charath\'eodory in 
the Calculus of Variations of Multiple Integrals}", in
Lect. Notes in Pure and Appl. Math. {\bf 100} (1985) 455-469

\bibitem{Ru2} 
H. Rund,
``{\it Integral Formul\ae Associated with the Euler-Lagrange Operator
of Multiple Integral Problems in the Calculus of Variation}",
\AE quationes Math. {\bf 11} (1974) 212-229

\bibitem{Sa} 
D. J. Saunders, 
``{\it An Alternative Approach to the Cartan Form in the Lagrangian Field 
Theories}", 
J. Phys. {\bf A 20} (1987) 339-349

\bibitem{Su} 
J. M. Souriau, 
``{\it Structure des Systemes Dynamique}", 
Dunod, Paris, 1970.

\end{thebibliography}
\end{document}